\documentclass[final,3p,times]{elsarticle}

\usepackage{epstopdf}
\usepackage{amssymb}
\usepackage{color}
\usepackage{mathtools}
\usepackage{xcolor,colortbl}
\definecolor{gray}{rgb}{0.82,0.82,0.82}

\usepackage{multirow}
\usepackage{multicol}

\newtheorem{theorem}{Theorem}[section]
\newtheorem{corollary}{Corollary}[theorem]
\newtheorem{lemma}[theorem]{Lemma}
\newtheorem{definition}[theorem]{Definition}

\newcommand{\RevA}[1]{{\color{black}{#1}}} 
\newcommand{\RevB}[1]{{\color{black}{#1}}} 
\newcommand{\RevUs}[1]{{\color{black}{#1}}}

\journal{Journal of Computational Physics}

\begin{document}

\begin{frontmatter}


\title{Time-Accurate and highly-Stable Explicit operators\\ for stiff differential equations}

\author{Maxime Bassenne\footnote{Corresponding Author. E-mail address: bassenne@stanford.edu. Now at Laboratory of Artificial Intelligence in Medicine and Biomedical Physics, Stanford University, Stanford, CA 94305.}, Lin Fu, Ali Mani}
\address{Center for Turbulence Research, Stanford University, Stanford, CA 94305-3024}%

\begin{abstract}
Unconditionally stable implicit time-marching methods are powerful in solving stiff differential equations efficiently. In this work, a novel framework to handle stiff physical terms implicitly is proposed. Both physical and numerical stiffness originating from convection, diffusion and source terms (typically related to reaction) can be handled by a set of predefined Time-Accurate and highly-Stable Explicit (TASE) operators in a unified framework. The proposed TASE operators act as preconditioners on the stiff terms and can be deployed to any existing explicit time-marching methods straightforwardly. The resulting time integration methods remain the original explicit time-marching schemes, yet with nearly unconditional stability. The TASE operators can be designed to be arbitrarily high-order accurate with Richardson extrapolation such that the accuracy order of original explicit time-marching method is preserved. Theoretical analyses and stability diagrams show that the $s$-stages $s$th-order explicit Runge-Kutta (RK) methods are unconditionally stable when preconditioned by the TASE operators with order $p \leq s$ and $p \leq 2$. On the other hand, the $s$th-order RK methods preconditioned by the TASE operators with order of $p \leq s$ and $p > 2$ are nearly unconditionally stable. The only free parameter in TASE operators can be determined a priori based on stability arguments. \RevUs{Unlike classical implicit methods, the TASE methodology allows for solving non-linear problems with arbitrary order without requiring solving a nonlinear system of equations.} A set of benchmark problems with strong stiffness is simulated to assess the performance of the TASE method. Numerical results suggest that the proposed framework preserves the high-order accuracy of the explicit time-marching methods with very-large time steps for all the considered cases. As an alternative to established implicit strategies, TASE method is promising for the efficient computation of stiff physical problems.

\end{abstract}

\end{frontmatter}

\section{Introduction} \label{sec1}

Many scientific and engineering problems can be modeled as partial differential equations (PDEs) or ordinary differential equations (ODEs), and solved numerically with the assistance of modern computations \cite{dafermos2000hyperbolic,anderson1995computational,kim1987turbulence}. The success of solving these problems strongly depends on the accurate, efficient, and stable numerical algorithms for both temporal and spatial discretizations. These are however not trivial due to the possible presence of stiffness originating from both numerics and physics \cite{curtiss1952integration,lomax2013fundamentals,mortazavi2015computational}. Numerical stiffness typically induced by regions of excessive spatial resolutions imposes that the time-step should be chosen even smaller than the smallest physical time scale
of the system. A common situation is when structured meshes are used in cylindrical coordinates leading to excessive azimuthal resolution near the axis of symmetry. Physical stiffness occurs when there is a separation of scales between a fast (quasi-)steady physical phenomenon and slower unsteady physical processes with larger time scales. In both cases, the common feature is that the numerical solution exhibits fast decaying modes that the user is not interested in.

In wide range of multi-physics systems the problematic stiff terms are not known ahead of designing the numerical algorithm, e.g. due to nonlinearities in the governing differential equations, or because of the lack of prior insights about the phenomena to be simulated such as the identification of competing physics. In such cases it would be convenient to first develop a fully explicit numerical solver since explicit codes are often simple to program and debug. In other cases, one may have quickly modified an already available explicit code to simulate a new problem that is later identified to be stiff. The stiff processes frequently originate from only a few terms in the governing equations that compete in a quasi-steady manner. These terms can be identified by examining the impact of each term on the behavior of the solution and their characteristic time scales. In these scenarios users often desire to insert a quick fix to handle these specific stiff terms, while maintaining the structure of an existing verified code, rather than re-programing the entire code using traditional implicit schemes. TASE operators are particularly useful as quick programming remedies for such scenarios, as in many of these conditions the programming time may overwhelm the simulation time itself.

For temporal discretization for stiff PDEs and ODEs, unconditionally stable implicit time-marching methods are currently a common choice as they enable a flexible choice of the time-step size based on the desired accuracy \cite{dahlquist1963special,ascher1995implicit,boom2018optimization}. Taking aerodynamic applications for instance, the quantities of interests, e.g. drag and lift, are only relevant to the steady solution of Euler or Navier-Stokes (NS) equations \cite{sclafani2008cfl3d,tinoco2017summary,rumsey2018overview}. For high computational efficiency, the nonlinear convection and diffusion terms of the governing laws can be handled implicitly, linearized in time with first-order approximation, and marched with the efficient Lower-upper symmetric-Gauss-Seidel (LU-SGS) \cite{yoon1988lower,chen2000fast,zhang2004block} method, which eliminates the need for block-diagonal inversions without using a diagonalization procedure. This first-order accurate implicit framework allows for stable computations of flows across a wide range of Mach numbers without time-step constraints. To recover the temporal accuracy for unsteady flows, a second-order implicit time-marching method, which is both A- and L-stable \cite{dahlquist1963special}, can be developed by further combing the dual time-stepping strategy \cite{jameson1991time,jameson2015application} with the second-order backward difference formula (BDF2). Another common numerical challenge is to accurately simulate the time dynamics of reactive flows with stiff chemistry, for which the time scale is determined by the physico-chemical process related to reaction and transport. In cases where the smallest chemical time scale is comparable to that of transport process, the coupled system can be marched simultaneously via explicit methods. In general however, the chemical time scale is orders of magnitude smaller than the transport-related time scale, therefore an explicit treatment of chemical terms leads to prohibitively small computational time steps and unaffordable computational costs. A popular strategy for this type of convection-diffusion-reaction equations is to treat the non-stiff convection-diffusion operators explicitly and the reaction-related stiff source terms implicitly by a class of operator-splitting methods \cite{marchuk1968some,strang1968construction}. The resulting ODEs for the reaction steps can be integrated efficiently by fully implicit multi-step variable-order BDF methods \cite{gear1971automatic} or Rosenbrock-Krylov methods \cite{tranquilli2014rosenbrock}. Although these splitting methods are strongly stable \cite{knio1999semi,singer2006operator}, the overall temporal accuracy is typically limited to second order \cite{speth2013balanced}.

As it is well known, for scientific computations with temporal accuracy requirements, higher-order methods allow for reaching a target accuracy at a lower overall computational cost, despite a larger computational cost per time step. Consequently, many efforts have been devoted to the development of unconditionally stable time-marching methods of high order (higher than second order), e.g. the cyclic and composite linear multistep methods \cite{donelson1971cyclic,sloate1973stable}, the advanced-step-point methods \cite{psihoyios1998stability,aiken1985stiff,cash1983integration}, and hybrid methods \cite{england1982some}. However, A-stable linear multi-step schemes are at most second-order accurate \cite{dahlquist1963special}. On the contrary, it is possible to derive both A- and L-stable implicit RK methods that yield higher-order accuracy with strong stability \cite{butcher1987numerical,butcher2003numerical,butcher1964implicit}. Implicit RK methods can be classified based on the form of their coefficient matrix (also called ``Butcher tableaux''). Specifically, the implicit RK methods characterized by a lower triangular coefficient matrix with nonzero diagonal elements are called diagonally-implicit RK method (DIRK) \cite{alt1971methodes,crouzeix1975approximation,alexander1977diagonally}. Compared to those with full coefficient matrix, the solution at each RK substep can be solved sequentially with high efficiency as it is independent of subsequent stage values \cite{alexander1977diagonally}. The efficiency can be further improved by the singly-diagonally implicit RK (\RevUs{SDIRK}) methods \cite{wanner1996solving}, for which the diagonal coefficients of the lower triangular matrix are identical and the operations on the Jacobian can be reused for each RK substep throughout these iterations without any significant reduction in convergence when the Newton's method is employed to solve the resulting nonlinear system \cite{crouzeix1975approximation}. 

As another class of efficient time-marching methods for convection-diffusion equations, the implicit-explicit (IMEX) \cite{ascher1997implicit,calvo2001linearly,ascher1995implicit,cooper1983additive,wang2015stability} methods consist of employing an implicit discretization for the diffusion terms and an explicit discretization for the convection terms. The rationale is that the diffusion terms are typically stiff and linear while the convection terms are not stiff but nonlinear. It is quite challenging to handle the nonlinear convection terms implicitly as the resulting nonlinear system needs to be solved by an iterative solver. When combined with the high-order linear multistep schemes \cite{donelson1971cyclic,sloate1973stable}, the IMEX methods still suffer from strong time-step constraint \cite{ascher1997implicit,calvo2001linearly}. Instead, the IMEX methods are typically developed based on the high-order RK method with larger stability regimes \cite{ascher1997implicit,wang2015stability}. An $s$-stages high-order IMEX-RK method \cite{ascher1997implicit,calvo2001linearly,wang2015stability} can be defined by deploying an $s$-stages DIRK \cite{alexander1977diagonally} method for diffusion terms and an $(s+1)$-stages explicit RK method for convection terms. \RevB{One main drawback of IMEX-RK methods is that the accuracy-order is typically limited to be lower than fifth due to the prohibitive difficulty in satisfying the higher-accuracy-order conditions when determining the coefficients. Alternative high-order IMEX methods have been proposed in \citep{constantinescu2010extrapolated} to handle problems with both stiff and nonstiff terms efficiently based on the extrapolated method, which allows for any desired accuracy order for ODEs \cite{deuflhard1985recent} and PDEs \cite{hairer1988extrapolation}. The uniqueness is that the extrapolated method is fairly friendly to parallelization and thus can enhance the efficiency on multicore hardware architectures, which is non-trivial for other implicit methods. The present work builds on similar extrapolation concepts to achieve higher-order accuracy, yet it is formulated at the semi-discrete level independent of a specific time-marching scheme, yielding provable stability properties.} IMEX time-marching methods have been widely used to solve different problems \cite{zhong1996additive,pareschi2005implicit,gottlieb2012stability,kennedy2003additive,wang2016local,fu2017analysis}.

In the present work, we propose an alternative framework to implicitly handle stiff ODEs or PDEs. A family of unique TASE operators with arbitrarily high-order accuracy is defined that can be deployed to any existing explicit time-marching methods. The resulting time integration remains explicit, yet it becomes nearly unconditional stable. The stiffness originating from different mechanisms, including convection, diffusion, and source terms, can be treated implicitly in the unified TASE framework. Moreover, the order of accuracy of the baseline explicit time-marching methods can be maintained. The stability properties of the TASE operators in combination with established RK methods is analyzed by theoretical proof and stability diagrams. The performance of TASE operators in terms of accuracy and stability are demonstrated by a set of canonical benchmark simulations.

The rest of the paper is organized as follows. The TASE framework is introduced in Section~\ref{TASE_operators},which also includes the analytical description of a family of operators. The stability properties of the proposed TASE operators are analyzed in Section~\ref{stability}. Section~\ref{Numerical_validations} further discuss the practical deployment of TASE operators to challenging physical problems with source terms, multiple stiff mechanisms, non-linear equations, and grid-induced stiffness in polar coordinates. Concluding remarks and future work are discussed in the last section.

\section{TASE operators} \label{TASE_operators}

The problem of formulating TASE operators is described in Section~\ref{sec_problem_definition}. A specific example of a first-order solution is described in Section~\ref{sec_first_TASE}, and a family of generalized TASE operators is introduced in Section~\ref{sec_high_TASE}. The free parameter, the explicit definitions, and the stability diagrams of TASE operators are discussed in Sections~\ref{sec_alpha}--\ref{sec_stability_digram}.

\subsection{Problem definition} \label{sec_problem_definition}

Typical theoretical descriptions of physical and engineering problems involve conservation equations represented as
\begin{equation}
\frac{\partial y}{\partial t}=f\left( x, t, y(x,t) \right), \label{eq1_1}
\end{equation}
where $(x, t)$ denote space and time, respectively, $y(x, t)$ is a vector of unknown physical variables, and $f$ is a function that describes the time-rate of change of $y$. The numerical solution of Eq.~(\ref{eq1_1}) is generally sought on a discrete grid $(X,T)=\{x_{i}, t_j\}$, whose resolution ought to be chosen carefully for the numerical solution to faithfully approximate the exact solution. This is a major requirement of any physics-based computational solver. The spatial mesh resolution $\Delta x$ should resolve all the physical length scales of the problem. These length scales are determined by the the problem geometry, the initial condition, and the physical features that may arise over time. After spatial discretization the continuous conceptual model is turned into a semi-discrete computer model
\begin{equation}
\frac{d Y}{d t}=F\left(X, t, Y(t) \right), \label{eq1_2}
\end{equation}
where the numerical approximation of the operator $f$ is denoted as $F$, and $Y$ represents a vector of spatially discretized unknowns. In order to obtain the numerical solution to Eq.~(\ref{eq1_2}) at different time, the time derivative is numerically approximated with a certain time-marching scheme and time step $\Delta t$. The temporal resolution is constrained by both physical and numerical requirements.

Consider the linear model problem for the semi-discrete ordinary differential equation given in Eq.~(\ref{eq1_2})
\begin{equation}
\frac{d Y}{d t} = L Y, \label{eq21_1}
\end{equation} 
where $Y(t)$ is a vector of unknowns, and $L$ is a matrix operator. The explicit time integration of Eq.~(\ref{eq21_1}) is subject to an upper-bound constraint for the maximum time step that can be used for the numerical solution to be stable, as described in Section~\ref{sec1}. The objective of the present paper is to propose a family of operators $T^{(p)}_L$ such that the explicit time integration of 
\begin{equation}
\frac{d Y}{d t} = \left( T^{(p)}_L L \right) Y \label{eq21_2}
\end{equation} 
is unconditionally or nearly unconditionally stable irrespective of the choice of the time-marching scheme, and requiring that
\begin{equation}
T^{(p)}_L = 1 + \mathcal{O}(\Delta t^p), \label{eq21_3}
\end{equation}
which ensures that the exact solution to Eq.~(\ref{eq21_2}) converges to the exact solution of Eq.~(\ref{eq21_1}) with integer order $p$ in the limit of small time steps. In other words, the goal is to enforce that the numerical solution obtained with the TASE framework collapses to that that would be obtained with a unmodified explicit scheme in the limit of small time-steps, consistently with the fact that it is likely that the explicit time integration is stable on its own in this limit. Additionally, when the order $p$ is chosen to match the order of accuracy of the time-stepping scheme, the global order of accuracy of the numerical solver is not affected by replacing $L$ in the code with $(T_L^{(p)} L)$ as in Eq.~(\ref{eq21_2}). 

\RevB{For linear problems,} the TASE operator $T_L^{(p)}$ can be viewed as a preconditioner that transforms the exact operator $L$ into a modified operator 
\begin{equation}
\mathcal{T}_L^{(p)} = T_L^{(p)} L, \label{eq21_4}
\end{equation}
such that $d Y/dt = \mathcal{T}_L^{(p)} Y$. For various reasons highlighted in the following sections, the analyses of the TASE framework can be more clearly clarified using $T_L^{(p)}$, therefore for the rest of \RevB{this section} we mostly retain $T_L^{(p)}$ in our notation. \RevB{We provide general explicit definitions of both TASE preconditioners $\mathcal{T}_L^{(p)}$ and operators $T_L^{(p)}$ in Section~\ref{sec_formulas}.} 

\RevB{For certain classes of problems such as non-autonomous equations with linear operator $L$, the proposed method remains applicable even though $L$ is time-dependent. Without loss of generality, we have omitted the time-dependence of $L$ in the rest of the paper to facilitate the presentation. Moreover, the above TASE framework is also applicable to nonlinear equations as it relies on the linearized operator without loss of accuracy order, as shown in Sections~\ref{stability} and \ref{nonlinear_problem}.}

\subsection{Example of a first-order TASE operator} \label{sec_first_TASE}

The time-marching of Eq.~(\ref{eq21_1}), $dY/dt = LY$, with implicit Euler is unconditionally stable, and yields the fully-discrete computer model
\begin{equation}
\frac{Y^{n+1}-Y^{n}}{\Delta t} = L Y^{n+1}, \label{eq22_1}
\end{equation} 
where $Y^n=Y(t_n)$ and $t_{n+1}-t_n=\Delta t$. Using the matrix equality $((1-\Delta t L)^{-1}-1)/\Delta t = (1-\Delta t L)^{-1} L$, and assuming that $(1-\Delta t L)$ is invertible, Eq.~(\ref{eq22_1}) can be re-arranged as
\begin{equation}
\frac{Y^{n+1}-Y^{n}}{\Delta t} = \left((1-\Delta t L)^{-1}  L\right)\ Y^n. \label{eq22_2}
\end{equation} 

Equation~(\ref{eq22_2}) resembles the fully-discrete equation that would be obtained by time-advancing Eq.~(\ref{eq21_2}), $dY/dt = (T^{(p)}_L L)Y$ with explicit Euler, yet with unconditionally stable property by construction. Term by term identification of both equations suggest an example of first-order solution to the TASE problem definition formulated in Section~\ref{sec_problem_definition}, i.e.
\begin{equation}
\RevA{T^{(1)}_L} = \left( 1 - \Delta t L \right)^{-1}. \label{eq22_3}
\end{equation}
The time-marching of Eq.~(\ref{eq21_2}), $dY/dt = (T^{(p)}_L L)Y$, using the TASE operator defined in Eq.~(\ref{eq22_3}) is unconditionally stable if the time-marching scheme is explicit Euler, by construction. Furthermore, it is trivially first-order since its Laurent series is $(1 - \Delta t L)^{-1}=\sum_{k=0}^{\infty} (\Delta t L)^k$ in the limit of small time step $\Delta t$ \cite{markushevich1977theory}. This first-order operator satisfies only partly the design constraints described in Section~\ref{sec_problem_definition}, because its stability properties for other time-marching schemes are still unknown, and it is only first-order. The extension of this operator to general explicit time-marching schemes of arbitrary accuracy order is described in the next section.

\subsection{Generalization and extension to arbitrary order} \label{sec_high_TASE}

\begin{table}[t]
\begin{center}
\begin{tabular}{| c | c | c | c | c |} 
  \cline{1-5}
   & \multicolumn{4}{c|}{\textbf{TASE order}} \\ 
  \cline{2-5}
   & $p=1$ & $p=2$ & $p=3$ & $p=4$ \\ 
  \hline
  RK1 $(C=2.00)$ & $\mathbf{0.50}$ & & & \\
  RK2 $(C=2.00)$ & 0.50 & $\mathbf{1.50}$ & & \\
  RK3 $(C=2.50)$ & 0.40 & 1.20 & $\mathbf{2.80}$ & \\
  RK4 $(C=2.79)$ & 0.36 & 1.08 & 2.51 & $\mathbf{5.38}$ \\ 
  \hline
\end{tabular}
\end{center}
\caption{Values of $\alpha_\mathrm{min}$ for common explicit RK methods preconditioned by the TASE operators of varying order. In each column, the bold value is the value that corresponds to matching the order of the TASE operator with that of the explicit RK scheme. Note that practically it is not useful to derive the tailored TASE operators with accuracy order higher than that of the original explicit time-marching methods.} \label{tab_alpha}
\end{table}

The TASE operator of order $p$ is recursively defined as
\begin{equation}
   T^{(p)}_L \RevA{(\alpha, \Delta t)}= \left \{
  \begin{aligned}
    & (1- \alpha \Delta t L)^{-1}, && \text{if}\ p=1 \\ 
    &  \frac{2^{p-1} T^{(p-1)}_L \RevA{(\alpha/2, \Delta t)} - T^{(p-1)}_L \RevA{(\alpha, \Delta t)}}{2^{p-1}-1}, && \text{if}\ p\geq2
  \end{aligned} \right. \label{eq11}
\end{equation} 
where $\alpha >0$ is a free parameter that depends on the time-marching scheme and must satisfy
\begin{equation}
\alpha\geq \alpha_{\mathrm{min}}= \frac{2^p-1}{C}. \label{eq23_1}
\end{equation}
The denominator $C$ is the maximum value of $|\lambda \Delta t|$ that guarantees a stable numerical solution for a given explicit time-marching scheme when applied to the model problem $dy/dt=\lambda y$, with $\lambda<0$ and real unknown $y$ \cite{eberly2008stability,lambert1991numerical}. It can also be interpreted as the absolute value of the intersection of the stability diagram of a given explicit time-marching scheme with the negative real axis, which is a often known value for common schemes. Depending on the time-marching scheme and desired TASE accuracy order, the order-unity minimum value of $\alpha$ can be easily obtained from Eq.~(\ref{eq23_1}). For instance, the values of $\alpha_\mathrm{min}$ for common explicit RK methods preconditioned by the TASE operators of varying order are given in Table~\ref{tab_alpha}. The motivation for the constraint given by Eq.~(\ref{eq23_1}) on the choice of $\alpha$ becomes clearer when analyzing the stability properties of the defined TASE operators in Section~\ref{stability}.

The specific first-order operator derived in Section~\ref{sec_first_TASE} corresponds to $T_{L}^{(1)}\RevA{(\alpha,\Delta t)}$ with $\alpha=1$ that is greater than the threshold value $\alpha_{\mathrm{min}}=0.5$ given in Table~\ref{tab_alpha} for the combination of the explicit Euler scheme with the first-order TASE operator. However, when deployed with a high-order time-marching scheme, it would reduce the global order of accuracy of the solver to one. For particular scenarios, where the steady solution is pursued, this degeneration of accuracy order is fine. However, for most common cases with temporal accuracy requirement, it is yet desirable that the accuracy order of explicit time-marching methods is preserved after preconditioned with the TASE operator. The family of TASE operators defined in Eq.~(\ref{eq11}) generalizes the example given in Section~\ref{sec_first_TASE} to arbitrarily high order via Richardson extrapolation \RevB{\cite{moin2010fundamentals}}, and a broader class of time-marching schemes via the parameter $\alpha$. \RevB{The linear combination of two $(p-1)$-th order TASE operators with appropriate $\alpha$-values leads to cancellation of the $(p-1)$-th order error terms, hence resulting in a $p$-th order leading error term.}

The proof of the statement of accuracy for the TASE operator given in Eq.~(\ref{eq21_3}) is straightforward due to the recursive application of Richardson extrapolation in Eq.~(\ref{eq11}) \cite{richardson1911ix,1927_Richardson}. Note that Richardson extrapolation has been previously employed to improve the accuracy order of the strongly stable $\theta$-method  to solve stiff ODEs implicitly \cite{2010_Zlatev_theta}. One secondary yet useful outcome of examining the Taylor series of Eq.~(\ref{eq21_3}) and Eq.~(\ref{eq11}) is the conclusion that the leading error term increases with $\alpha$, which indicates that $\alpha=\alpha_{\mathrm{min}}$ is a good choice in practice.

\subsection{Rationale for the minimum value of $\alpha$} \label{sec_alpha}
A necessary condition for the time integration of $dy/dt = T^{(p)}_{\lambda}\RevA{(\alpha, \Delta t)}\lambda y$ to be unconditionally stable is that the modified eigenvalue $T^{(p)}_{\lambda}\RevA{(\alpha, \Delta t)}\lambda$ falls within the stability diagram of the given explicit time-marching scheme as the time step goes to infinity. We show below that this necessary condition collapses to the constraint prescribed on $\alpha$ in Eq.~(\ref{eq23_1}),
\begin{theorem} [Asymptotic stability]
\label{asymptotic_TASE}
The arbitrary order explicit time integration of $dy/dt = T^{(p)}_{\lambda}\RevA{(\alpha, \Delta t)}\lambda y$ remains stable in the limit of very large time steps when $\alpha\geq \alpha_{\mathrm{min}}= (2^p-1)/C$ (necessary condition for unconditional stability).
\end{theorem}

It follows from the Laurent series of $T^{(1)}_\lambda$ in the limit of large time steps $\Delta t$ and recursive application of Eq.~(\ref{eq11}), that 
\begin{equation}
\left( T^{(p)}_\lambda\RevA{(\alpha, \Delta t)} \lambda \right) \Delta t  = -\frac{2^p-1}{\alpha} + \mathcal{O}(\Delta t^{-1}). \label{eq24_1}
\end{equation} 
Using the definition of $\alpha_{\mathrm{min}}$ given in Eq.~(\ref{eq23_1}), the right-hand side term becomes $-(2^p-1)/\alpha = -(\alpha_{\mathrm{min}}/\alpha) C$. In conjunction with Eq.~(\ref{eq24_1}),
it implies that in the limit of large $\Delta t$,
\begin{equation}
- C \leq \left( T^{(p)}_\lambda\RevA{(\alpha, \Delta t)} \lambda \right) \Delta t < 0. \label{eq24_2}
\end{equation}
Equation~(\ref{eq24_2}) proves that in the limit of large time steps $\Delta t$, the modified dimensionless eigenvalue $\left( T^{(p)}_\lambda\RevA{(\alpha, \Delta t)} \lambda \right) \Delta t$ falls on the real axis and within the stability diagram of all time-marching schemes whose stability diagram cross the real axis at $(-C,0)$. This ensures stability in the asymptotic limit of large $\Delta t$, which concludes the proof of theorem (\ref{asymptotic_TASE}).  The above asymptotic analysis shows that the prescribed real value of $\alpha$ given in Eq.~(\ref{eq23_1}) is not compatible with time-marching schemes whose stability diagram is restricted to purely imaginary $\lambda$. These schemes are excluded from the discussion in the present paper anyway, since they are only useful for purely oscillatory physics which require to be time resolved and consequently rule out the need for their implicit treatment.
In Section~\ref{stability}, we will show that the necessary condition given in the present section for $\alpha$ works well in practice.

\subsection{Explicit definition of TASE operators} \label{sec_formulas}

The formulation of TASE operator only includes the parameter $\alpha$, the time step $\Delta t$ and the operator $L$. The explicit expressions can be obtained by \RevB{expanding} the recurrence formula prescribed in Eq.~(\ref{eq11}). \RevA{Although the application of TASE operators involves matrix inversions, which are reminiscent of implicit \textit{time-advancement} schemes, TASE \textit{operators} belong to the class of explicit functions as they are fully defined based on known values of the problem variables.} \RevB{For practical purposes, the TASE preconditioners $T^{(p)}_L$ and operators $\mathcal{T}^{(p)}_L$ are explicitly defined up to order $p=4$ below. The best suited formulation depends on the nature of the problem and practical implementation details. For linear problems, the operator form $\mathcal{T}^{(p)}_L$ may be better suited as computing $T_L^{(p)}L Y$ requires an additional matrix-vector product besides solely solving $p$ systems of linear equations as in $\mathcal{T}_L^{(p)} Y$. For nonlinear problems, as detailed in Section~\ref{nonlinear_problem}, the method requires using the preconditioner form $T^{(p)}_L$ since $\mathcal{T}^{(p)}_L$ is not applicable.}

The TASE preconditioners $T^{(p)}_L$ are expressed as 
\begin{equation}
T_L^{(p)}\RevA{(\alpha, \Delta t)} = \sum_{k=0}^{p-1} \beta_{p,k} \left(2^{k}-\alpha \Delta t L\right)^{-1} \label{A1_eq1}
\end{equation}
where the coefficients $\beta_{p,k}$ are given in Table~\ref{A1_t1}.

\begin{table}[!ht]
\begin{center}
\begin{tabular}{ c | c | c | c | c} 
 & $\RevB{k=0}$ & $\RevB{k=1}$ & $\RevB{k=2}$ & $\RevB{k=3}$ \\ [0.5ex]
  \hline
  $\RevB{p=1}$ & 1 & - & - & - \\ 
  $\RevB{p=2}$ & -1 & 4 & - & - \\ 
  $\RevB{p=3}$ & 1/3 & -4 & 32/3 & - \\ 
  $\RevB{p=4}$ & -1/21 & 4/3 & -32/3 & 512/21 \\ 
\end{tabular}
\end{center}
\caption{\RevUs{Coefficients} $\beta_{p,k}$ in the definition of $T_L^{(p)}\RevA{(\alpha, \Delta t)}$ as prescribed by Eq.~(\ref{A1_eq1}).} \label{A1_t1}
\end{table}

The TASE operators $\mathcal{T}^{(p)}_L$ are expressed as 
\begin{equation}
\mathcal{T}_L^{(p)}\RevA{(\alpha, \Delta t)} = \frac{1}{\alpha \Delta t} \left( -(2^{p}-1) + \sum_{k=0}^{p-1} \gamma_{p,k} \left(2^{k}-\alpha \Delta t L\right)^{-1} \right) \label{A1_eq2}
\end{equation}
where the coefficients $\gamma_{p,k}$ are given in Table~\ref{A1_t2}.

\begin{table}[!ht]
\begin{center}
\begin{tabular}{ c | c | c | c | c} 
& $\RevB{k=0}$ & $\RevB{k=1}$ & $\RevB{k=2}$ & $\RevB{k=3}$ \\ [0.5ex]
  \hline
  $\RevB{p=1}$ & 1 & - & - & - \\ 
  $\RevB{p=2}$ & -1 & 8 & - & - \\ 
  $\RevB{p=3}$ & 1/3 & -8 & 128/3 & - \\ 
  $\RevB{p=4}$ & -1/21 & 8/3 & -128/3 & 4096/21 \\ 
\end{tabular}
\end{center}
\caption{\RevUs{Coefficients} $\gamma_{p,k}$ in the definition of $\mathcal{T}_L^{(p)}\RevA{(\alpha, \Delta t)}$ as prescribed by Eq.~(\ref{A1_eq2}).} \label{A1_t2}
\end{table}

In terms of the implementation of the TASE operator, $T_L^{(p)}$, or equivalently $\mathcal{T}_L^{(p)}$, does not need to be computed explicitly, rather it should be avoided since it is often a full matrix. The only required routine is the one that computes $\left( T_L^{(p)} L \right) Y$, or equivalently $\mathcal{T}_L^{(p)} Y$. For many practical operators $L$, this leads to a banded matrix to solve which is much less costly than directly computing the inverse of a matrix \cite{banerjee2014linear}. In general, the direct evaluation of $T_L^{(p)}$ requires computing the inverse of $p$ matrices, which makes the cost of the present method at least comparable to that of other established implicit schemes. Nevertheless, the computational savings allowed by the use of a larger time step can outweigh the additional cost induced by the TASE operators at each time step.

\subsection{Stability diagrams}\label{sec_stability_digram}
%

\begin{figure}[htbp]
    \vskip 0.1in
    \centering
    \includegraphics[trim={4.8cm 2cm 3.8cm 0cm}, width=\textwidth,angle=0,clip]{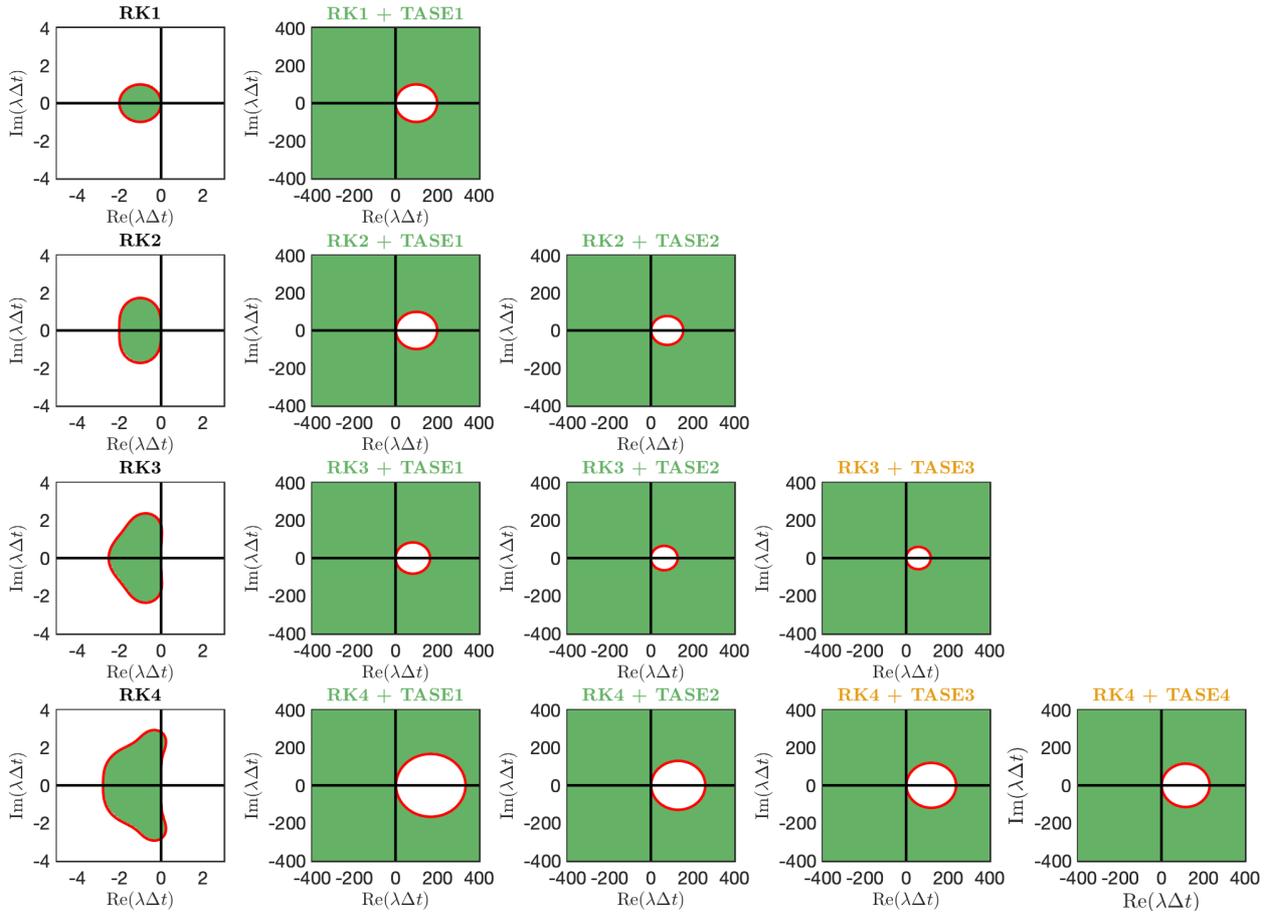}
    \caption{Stability diagrams of the explicit RK methods and the preconditioned RK methods with TASE operators up to fourth-order. The parameter $\alpha$ in TASE operators is set as the minimum value defined in Table.~\ref{tab_alpha}.}\label{Fig_RK_stability}
\vskip 0.1in
\end{figure}

\begin{figure}[htbp]
    \vskip 0.1in
    \centering
    \includegraphics[trim={4.1cm 2cm 3.5cm 0cm}, width=\textwidth,angle=0,clip]{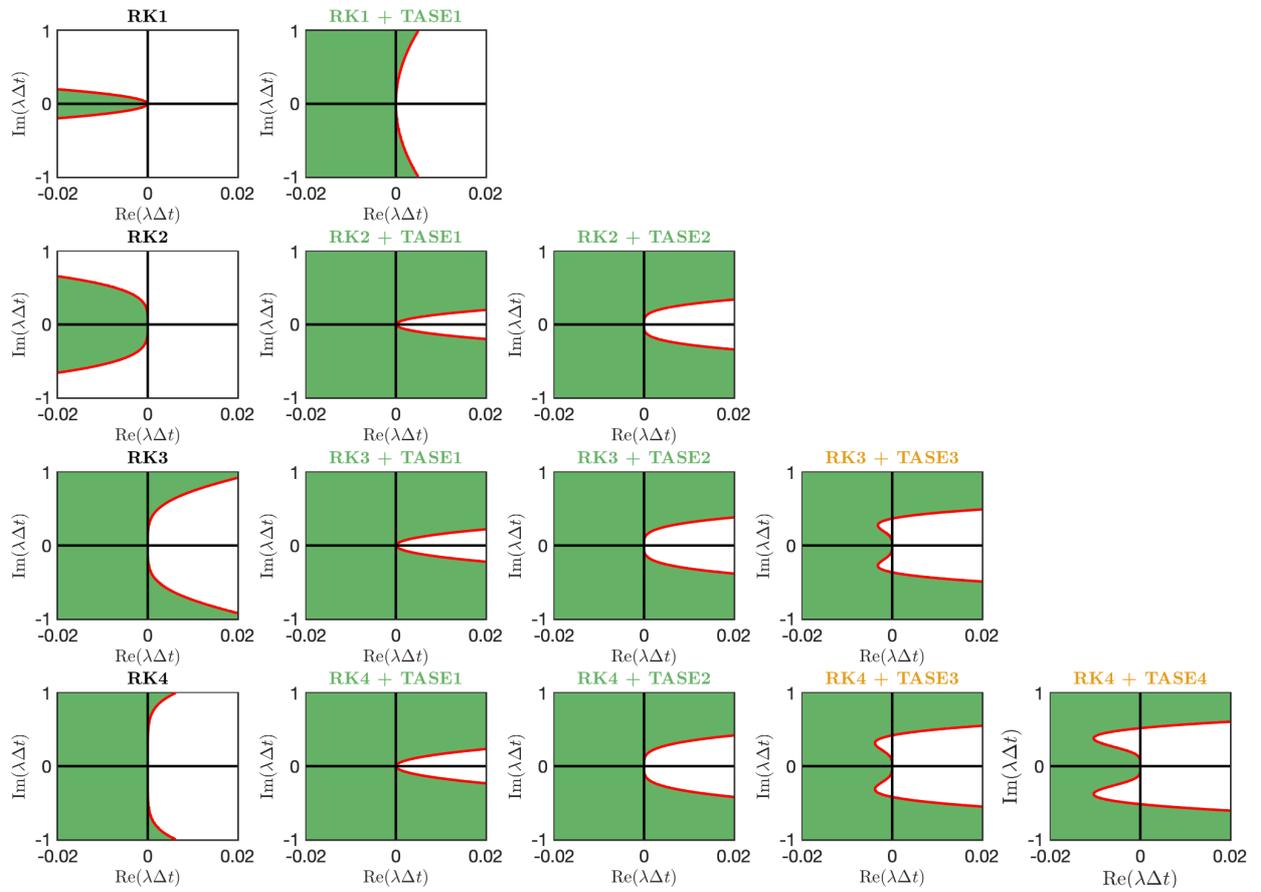}
    \caption{Zoomed-in view of the stability diagrams of the explicit RK methods and the preconditioned RK methods with TASE operators up to fourth-order. The parameter $\alpha$ in TASE operators is set as the minimum value defined in Table.~\ref{tab_alpha}.}\label{Fig_RK_stability_zoom}
\vskip 0.1in
\end{figure}

The stability diagrams for a set of $s$-stages $s$th-order RK methods preconditioned by TASE operators up to fourth-order are given in Fig.~\ref{Fig_RK_stability} and Fig.~\ref{Fig_RK_stability_zoom}. As shown in Fig.~\ref{Fig_RK_stability}, the stability regions of the explicit RK methods are substantially enlarged by the corresponding TASE operators when the parameter $\alpha$ is set as the minimum value defined in Table.~\ref{tab_alpha}. Specifically, all the RK methods preconditioned by the first- and second-order TASE operators are unconditionally stable, i.e. the entire left-half complex plane is mapped into the interior of the stability regions of the original explicit RK methods with TASE operators. On the other hand, the preconditioned RK methods become nearly unconditionally stable when combined with the TASE operators of order higher than second, as revealed in the zoomed-in Fig.~\ref{Fig_RK_stability_zoom}. Nevertheless, these unstable regions are tiny and correspond to operators strongly dominated by convection, as illustrated by the axes range in Fig.~\ref{Fig_RK_stability_zoom} that shows two orders of magnitude larger imaginary eigenvalue compared to their real counterparts. This class of problems is not of primary concern in the present analysis since stiffness originating from convection alone cannot be alleviated by the use of large time steps. Therefore, the nearly unconditionally stable RK methods may still allow for very-large time step in practice and will be demonstrated in the following numerical validations in Section~\ref{Numerical_validations}.

\begin{figure}[htbp]
    \vskip 0.1in
    \centering
    \includegraphics[trim={11.2cm 0cm 3.8cm 0cm}, width=\textwidth,angle=0,clip]{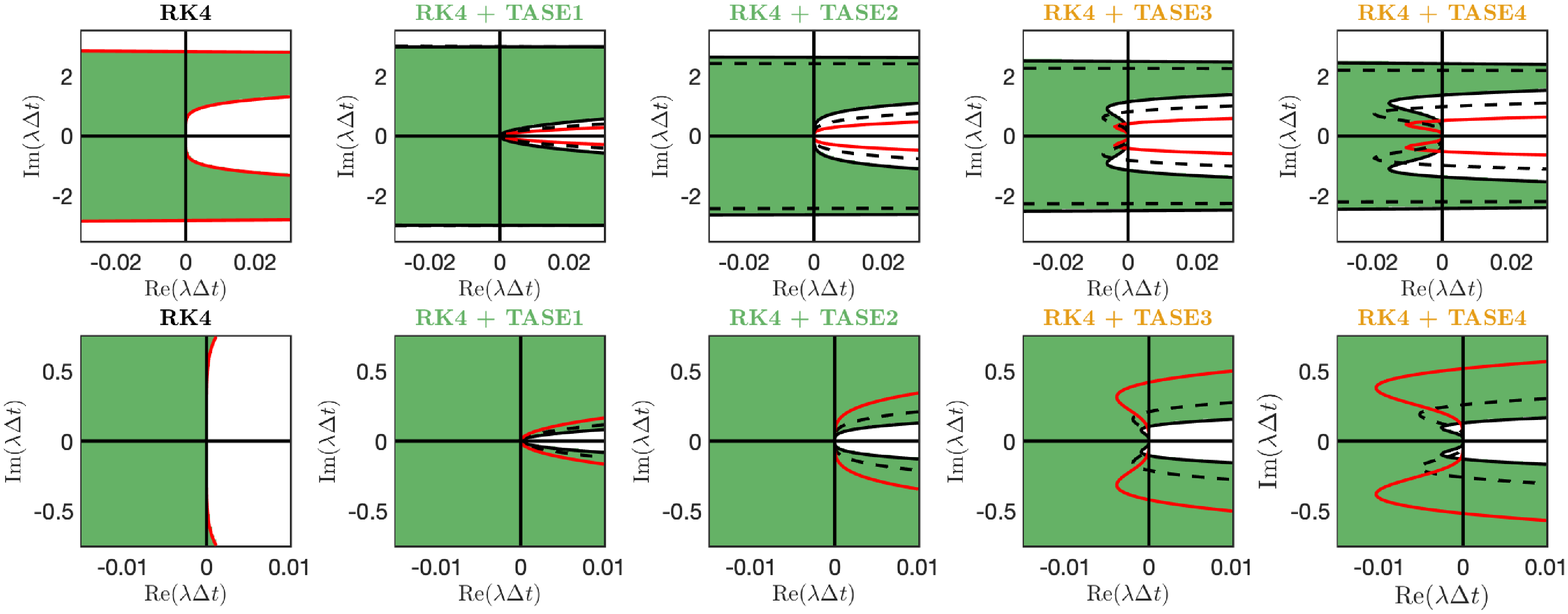}
    \caption{Zoomed-in views of the stability diagrams of the fourth-order RK method preconditioned by the TASE operators of first- to fourth-order. Top: stability boundaries from TASE operators with $\alpha =$ $0.25\alpha_\mathrm{min}$ (black solid lines) and $0.5\alpha_\mathrm{min}$ (black dashed lines); bottom: stability boundaries from TASE operators with $\alpha =$ $4\alpha_\mathrm{min}$ (black solid lines) and $2\alpha_\mathrm{min}$ (black dashed lines). The red lines for all plots denote the stability boundaries from the TASE operators with $\alpha =\alpha_\mathrm{min}$. The stability regions are colored by the green contours when $\alpha =$ $0.25\alpha_\mathrm{min}$ (top) and $4\alpha_\mathrm{min}$ (bottom) for TASE operators.}\label{Fig_RK_alpha}
\vskip 0.1in
\end{figure}

In order to investigate the dependence of the stability properties of the preconditioned RK methods on the choice of $\alpha$, the zoomed-in stability diagrams of the fourth-order RK method preconditioned by the TASE operators of first- to fourth-order with varying $\alpha$ are plotted in Fig.~\ref{Fig_RK_alpha}. It is shown that all the preconditioned RK methods have large unstable regions in the left-half complex plane when $\alpha < \alpha_\mathrm{min}$. Meanwhile, when $\alpha > \alpha_\mathrm{min}$, the stability of the preconditioned RK methods is improved with larger $\alpha$, i.e. the tiny unstable region if it exists in the left-half complex plane becomes even smaller.

\section{Analyses of stability properties} \label{stability}

In this section, the stability properties of the proposed TASE operators applied to several established explicit time-marching methods are analyzed in detail. The notion of stability chosen in this paper is the traditional notion of A-stability, which enforces that the numerical method be stable when the exact numerical solution is stable too, as defined below \cite{2010_Zlatev_theta,wanner1996solving}. \RevA{The A-stability property relies on the magnitude of the stability function $\sigma = y^{n+1}/y^n$, which is defined as the ratio of two successive numerical solution values \cite{moin2010fundamentals}.}
\begin{definition} [A-stability]
\label{Astability}
Consider the set $S$ containing all values of $z \in \mathbb{C}$ for which $\left| \sigma(z) \right| \leq 1$. A method with stability function $\sigma = y^{n+1}/y^n$ is said to be A-stable if $S \supset \mathbb{C}^- = \{z\ |\ \mathrm{Re}(z) \leq 0\}$.
\end{definition}

This section proves that the family of TASE operators defined in Section~\ref{sec_high_TASE} satisfies the constraint that the time-integration of Eq.~(\ref{eq21_2}), $dY/dt = \left( T^{(p)}_L L \right) Y$, is either unconditionally or nearly unconditionally stable. The analysis considers a single eigenvalue $\lambda$ of $L$, and its TASE substitute $\left( T^{(p)}_\lambda \lambda \right)$.  It is sufficient to do the analysis for individual eigenvalues of $L$ since the result extends to general matrices $L$ by eigenvalue decomposition, as the eigenvectors basis of modified operator $\left( T^{(p)}_L L \right)$ is the same as that of the original operator $L$.

In particular, we restrict the analysis to general $s$-stages $s$th-order RK methods, written in Shu-Osher form \cite{shu2002survey}:
\begin{equation}
\left \{
\begin{aligned}
y^{(0)} &= y^{n}, \\
y^{(i)} &= \sum_{k=0}^{i-1} \left( a_{i,k} y^{(k)} + \Delta t b_{i,k} f(y^{(k)}, t_n + c_k \Delta t) \right),\ \ \ a_{i,k} \geq 0,\ \ \ i=1, \dots, s \\
y^{n+1} &= y^{(s)}, \\
\end{aligned}
\right. \label{ShuRK}
\end{equation}
where by consistency $\sum_{k=0}^{i-1} a_{i,k}=1$. We additionally assume that all $b_{i,k}$'s are non-negative, i.e. $b_{i,k} \geq 0$, which is the case for most RK methods including $s$-stages $s$th-order with $s \leq 4$ ones of interest in the present study. For these specific methods, the stability diagram does not depend on the coefficients themselves, only on the order, or equivalently the number of stages, therefore there is no need to specify the exact coefficients here.

\subsection{Sufficient condition for unconditional stability of first-order TASE operators with arbitrary order RK} \label{proof_T1_RK}

We provide below the theoretical proof that the combination of first-order TASE operators with arbitrary order RK methods (\ref{ShuRK}) is unconditionally stable if $\alpha$ is chosen larger than a critical value.

\begin{theorem} [Unconditional stability of first-order TASE with RK]
\label{TASE1_RK}
The explicit RK time integration (\ref{ShuRK}) of $dy/dt = T^{(1)}_{\lambda}\RevA{(\alpha, \Delta t)}\lambda y$ is unconditionally stable when $\alpha \geq 0.5\  \underset{i,k}{\max}{(b_{i,k}/a_{i,k})}$ (sufficient condition).
\end{theorem}

We prove theorem \ref{TASE1_RK} in two steps. We first prove the particular case of explicit Euler time integration:
\begin{lemma}[First-order TASE: Explicit Euler]
\label{TASE1_RK1}
The explicit Euler time integration of $dy/dt = T^{(1)}_{\lambda}\RevA{(\alpha, \Delta t)}\lambda y$ is unconditionally stable, i.e. $\|y^{n} + \Delta t T^{(1)}_\lambda\RevA{(\alpha, \Delta t)}\lambda y^{n}\| \leq \|y^{n}\|$, if and only if $\alpha \geq 0.5$.
\end{lemma}

The explicit Euler time integration of $dy/dt = T^{(1)}_\lambda\RevA{(\alpha, \Delta t)} \lambda y$ yields $y^{n+1} = \left( 1 + \Delta t T_\lambda^{(1)}\RevA{(\alpha, \Delta t)}\lambda \right)y^n$. Using the definition of $T_\lambda^{(1)}\RevA{(\alpha, \Delta t)}$ in Eq.~(\ref{eq22_3}), it can be re-arranged as 
\begin{equation}
\frac{y^{n+1}-y^n}{\Delta t} = \alpha \lambda y^{n+1} + (1-\alpha) \lambda y^n. \label{eq25_2}
\end{equation}
The corresponding stability function $\sigma=y^{n+1}/y^n$ is \begin{equation}
\sigma = \frac{1+(1-\alpha)\lambda \Delta t}{1 - \alpha \lambda \Delta t}. \label{eq25_3}
\end{equation}
It is trivial to show that the condition $| \sigma | \leq 1$ subject to $\mathrm{Re}(\lambda) \leq 0$ is satisfied if and only if $\alpha \geq 0.5$. This concludes the proof of lemma \ref{TASE1_RK1}.

We now prove that unconditional stability for higher-order RK schemes directly follows from unconditional stability with explicit Euler:
\begin{lemma}[First-order TASE: From explicit Euler to arbitrary order RK]
\label{TASE1_RK_induction}
The arbitrary order RK time integration (\ref{ShuRK}) of $dy/dt = T^{(1)}_{\lambda}\RevA{(\alpha, \Delta t)}\lambda y$ is unconditionally stable if the explicit Euler integration of $dy/dt = T^{(1)}_{\lambda}\RevA{(\alpha_{i,k}, \Delta t)}\lambda y$ where $\alpha_{i,k} = (a_{i,k}/b_{i,k}) \alpha$ is unconditionally stable. 
\end{lemma}

The explicit time integration of $dy/dt = T^{(1)}_\lambda\RevA{(\alpha, \Delta t)} \lambda y$ with an arbitrary order RK method (\ref{ShuRK}) yields 
\begin{equation}
\left \{
\begin{aligned}
y^{(0)} &= y^{n}, \\
y^{(i)} &= \sum_{k=0}^{i-1} a_{i,k} \left(y^{(k)} + \Delta t T^{(1)}_\lambda\RevA{(\alpha, \Delta t)} \frac{b_{i,k}}{a_{i,k}} \lambda y^{(k)} \right),\ \ \ a_{i,k} \geq 0,\ \ \ i=1, \dots, s \\
y^{n+1} &= y^{(s)}, \\
\end{aligned}
\right. \label{ShuRK_linear}
\end{equation}

The key is to note from the definition of the first-order TASE operator in Eq.~(\ref{eq22_3}) that
\begin{equation}
T^{(1)}_\lambda\RevA{(\alpha, \Delta t)} \frac{b_{i,k}}{a_{i,k}} \lambda 
= T^{(1)}_{\lambda'}\RevA{(\alpha_{i,k}, \Delta t)} \lambda'
\end{equation}
where $\alpha_{i,k}=(a_{i,k}/b_{i,k})\alpha$ and $\lambda' = (b_{i,k}/a_{i,k}) \lambda$. Importantly, note that since $b_{i,k}/a_{i,k} \geq 0$, $\text{Re}(\lambda') \leq 0$ if and only if $\text{Re}(\lambda) \leq 0$. Therefore Eq.~(\ref{ShuRK_linear}) becomes
\begin{equation}
\left \{
\begin{aligned}
y^{(0)} &= y^{n}, \\
y^{(i)} &= \sum_{k=0}^{i-1} a_{i,k} \left(y^{(k)} + \Delta t T^{(1)}_{\lambda'}\RevA{(\alpha_{i,k}, \Delta t)} \lambda' y^{(k)} \right),\ \ \ \alpha_{i,k} \geq 0,\ \ \ i=1, \dots, s \\
y^{n+1} &= y^{(s)}, \\
\end{aligned}
\right. \label{ShuRK_linear_modif}
\end{equation}

It follows from Eq.~(\ref{ShuRK_linear_modif}) that all intermediate stages of arbitrary order RK schemes of the form (\ref{ShuRK}) are convex combinations of explicit Euler operations, with combination factors $a_{i,k}$ summing to 1, eigenvalue $\lambda$ replaced by $\lambda'$, and $\alpha$ replaced by $\alpha_{i,k}$. Therefore $\|y^{(0)}\| = \|y^n\|$ and $\|y^{(i)}\| \leq \sum_{k=0}^{i-1} a_{i,k} \left\| y^{(k)} + \Delta t T^{(1)}_{\lambda'}\RevA{(\alpha_{i,k}, \Delta t)} \lambda' y^{(k)} \right\|$ for $i=1,\dots, s$. If we assume that the explicit Euler integration of $dy/dt = T^{(1)}_{\lambda}\RevA{(\alpha_{i,k}, \Delta t)}\lambda y$ is unconditionally stable, it trivially follows that $\|y^{n+1}\| \leq \| y^n \|$, which concludes the proof of lemma \ref{TASE1_RK_induction}.

The complete proof of theorem \ref{TASE1_RK} is obtained by combining both lemmas. From lemma \ref{TASE1_RK_induction}, we conclude that for the explicit RK time integration (\ref{ShuRK}) of $dy/dt = T^{(1)}_{\lambda}\RevA{(\alpha, \Delta t)}\lambda y$ to be unconditionally stable, it is sufficient that the explicit Euler integration of $dy/dt = T^{(1)}_{\lambda}\RevA{(\alpha_{i,k}, \Delta t)}\lambda y$ where $\alpha_{i,k} = (a_{i,k}/b_{i,k}) \alpha$ is unconditionally stable. From lemma \ref{TASE1_RK}, we conclude that the explicit Euler integration of $dy/dt = T^{(1)}_{\lambda}\RevA{(\alpha_{i,k}, \Delta t)}\lambda y$ is unconditionally stable if only if $\alpha_{i,k} \geq 0.5$ for all $i$ and $k$, which is equivalent to $\alpha \geq 0.5 \max_{i,k}{(b_{i,k}/a_{i,k})}$. This concludes the proof of theorem \ref{TASE1_RK}.

In this section, we have proven a sufficient condition for the unconditional stability of the combination of  first-order TASE with arbitrary order RK. Given a specific RK method, one can write it in Shu-Osher form (\ref{ShuRK}) and compute the critical value $\alpha_{\mathrm{min}} = 0.5 \max_{i,k}{(b_{i,k}/a_{i,k})}$. For example, linear SSP RK methods up to order 4 \cite{gottlieb2001strong} are characterized by $\max_{i,k}{(b_{i,k}/a_{i,k})}=1$ hence $\alpha_{\mathrm{min}}=0.5$ irrespective of the order. Interestingly the resulting sufficient condition is more restrictive than the necessary condition derived in Section~\ref{sec_alpha} and shown in Table~\ref{tab_alpha} for RK3 and RK4. Stability diagrams shown in Section~\ref{sec_stability_digram}, and  numerical examples in Section~\ref{Numerical_validations} indicate that in practice the necessary condition provided in provided in Section~\ref{sec_high_TASE} and in Table~\ref{tab_alpha} works well. We further assess this in the proof below in Section~\ref{general_proof}, for first-order TASE but also for higher-order.

\subsection{General proof} \label{general_proof}

The analysis conducted in Sec.~\ref{proof_T1_RK} becomes cumbersome when the order of the TASE operator is larger than one. Here instead, we proceed with an alternative approach that relies on the following corollary \cite{2010_Zlatev_theta}
\begin{corollary}[A-stability proof from maximum modulus principle]
\label{Astability_proof}
A numerical method is A-stable if and only if\\
(a) $\sigma(z)$ is analytic in $\mathbb{C}^-$, and \\
(b) it is stable on the imaginary axis, i.e. $\left| \sigma(i y) \right| \leq 1$ for all real values of $y$.
\end{corollary}
The proof is made in two steps. First, we prove that $\sigma(z)$ is an analytic function in $\mathbb{C}^-$, then we assess the stability of the methods for purely imaginary eigenvalues.

(a) We first compute the stability function $\sigma(z)$ with $z=\lambda \Delta t$. The explicit time integration of $dy/dt=T^{(p)}_\lambda \RevA{(\alpha, \Delta t)} \lambda y$ with arbitrary order TASE and arbitrary order RK method yields
\begin{equation}
\left \{
\begin{aligned}
y^{(0)} &= y^{n}, \\
y^{(i)} &= \sum_{k=0}^{i-1} a_{i,k} \left(y^{(k)} + \Delta t T^{(p)}_\lambda\RevA{(\alpha, \Delta t)} \frac{b_{i,k}}{a_{i,k}} \lambda y^{(k)} \right),\ \ \ a_{i,k} \geq 0,\ \ \ i=1, \dots, s \\
y^{n+1} &= y^{(s)}, \\
\end{aligned}
\right. \label{ShuRK_linear_general}
\end{equation}
which can be recasted in the more explicit form
\begin{equation}
y^{n+1} = \left( 1 + \sum_{k=0}^{s-1} A_{s,k} \left(\Delta t T^{(p)}_\lambda\RevA{(\alpha, \Delta t)} \lambda\right)^{k+1} \right) y^n, \label{Shu_all}
\end{equation}
where $A_{s,k}$ can be recursively defined using linear combinations of the coefficients $a_{i,k}$ and $b_{i,k}$ \cite{gottlieb2001strong}. The exact definition of $A_{m,k}$ is irrelevant to the present proof, hence not provided here. Injecting the explicit definition of the TASE operators $T^{(p)}_\lambda$ given in Eq.~\ref{A1_eq1}, $T_\lambda^{(p)}\RevA{(\alpha, \Delta t)} = \sum_{i=0}^{p-1}\beta_{p,i}\left(2^{i}-\alpha \Delta t \lambda\right)^{-1}$ with coefficients $\beta_{p,k}$ given in Table~\ref{A1_t1}, we obtain the following stability function
\begin{equation}
\sigma(z) = \left( 1 + \sum_{k=0}^{s-1} A_{s,k} z^{k+1} \left( \sum_{i=0}^{p-1} \beta_{p,i} \left(2^{i}-\alpha z\right)^{-1} \right)^{k+1} \right). \label{injected}
\end{equation}
Using the multinomial theorem, Eq.~(\ref{injected}) can be rearranged as
\begin{equation}
\sigma(z) = \left( 1 + \sum_{k=0}^{s-1} A_{s,k} z^{k+1} \sum_{\ell_1+\ell_2+\dots+\ell_p=k+1} \binom{k+1}{\ell_1,\ell_2,\dots,\ell_p} \prod_{i=0}^{p-1} \beta_{p,i}^{\ell_i} \left(2^{i}-\alpha z \right)^{-\ell_i} \right). \label{injected_multinomial}
\end{equation}

It is possible to express the above stability function as $\sigma(z) = P(z) / Q(z)$ where $P$ and $Q$ are defined as 
\begin{equation}
Q(z) = \prod_{i=0}^{p-1} \left(2^{i}-\alpha z \right)^{s}, 
\end{equation}
and
\begin{equation}
P(z) = \left( Q(z) + \sum_{k=0}^{s-1} A_{s,k} z^{k+1} \sum_{l_1+l_2+\dots+l_p=k+1} \binom{k+1}{\ell_1,\ell_2,\dots,\ell_p} \prod_{i=0}^{p-1} \beta_{p,i}^{\ell_i} \left(2^{i}-\alpha z\right)^{s-\ell_i} \right). 
\end{equation}

The key here is to note that $s - \ell_i \geq 0$ hence $P(z)$ is a polynomial and, furthermore, \RevUs{the roots $z=2^i/\alpha \geq 0$ of $Q(z)$} are not in $\mathbb{C}^-$. Therefore $Q(z)$ is an analytic function in $\mathbb{C}^-$. As a polynomial, $P(z)$ is also an analytic function in $\mathbb{C}^-$. We can then conclude that as the ratio of two analytic functions where the denominator that does not admit poles in $\mathbb{C}^-$, the stability function \RevA{$\sigma(z)$} is proved to be an analytic function in $\mathbb{C}^-$ for arbitrary order TASE and arbitrary order RK scheme.

\begin{figure}[t]
    \vskip 0.1in
    \centering
    \includegraphics[width=\textwidth,angle=0]{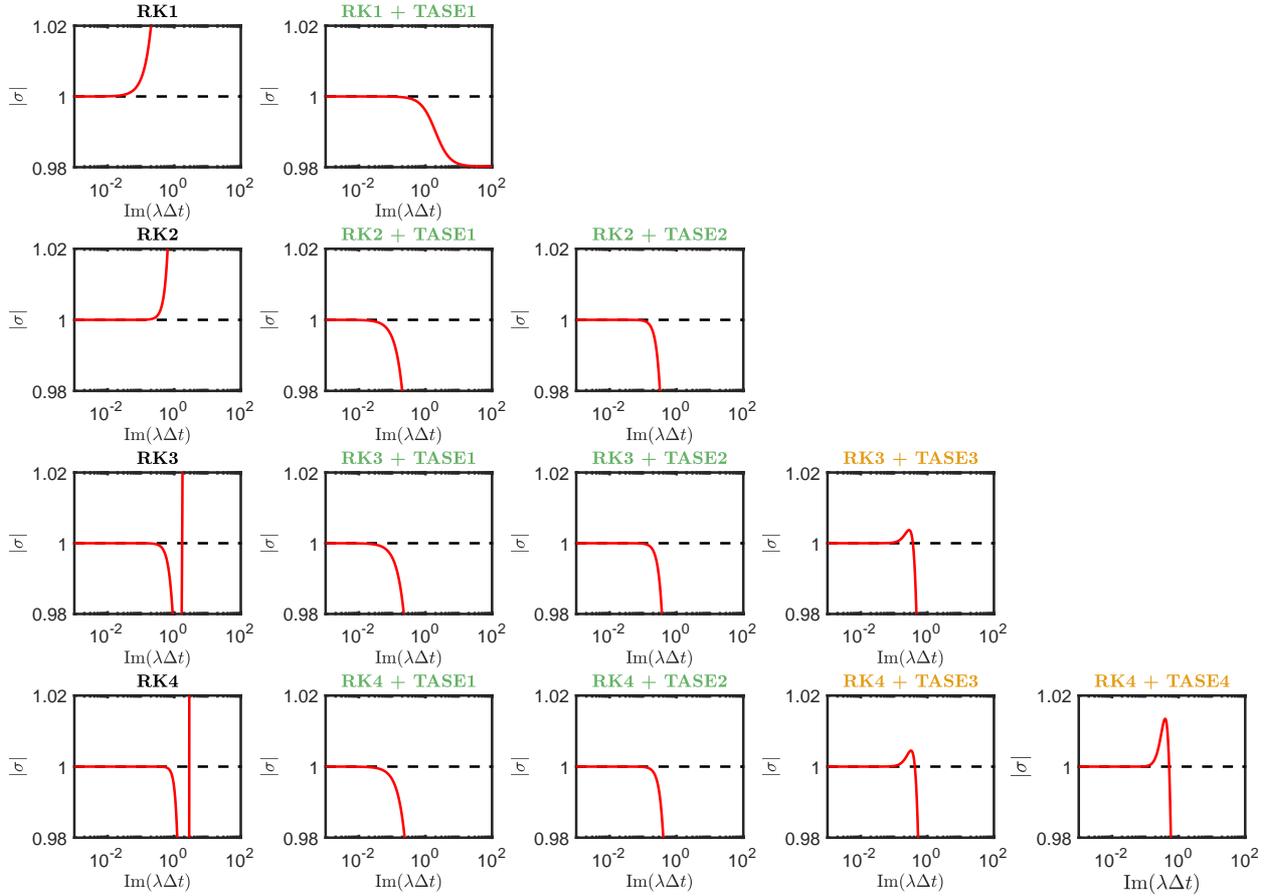}
    \caption{Zoomed-in view of the stability function values for purely imaginary eigenvalues $\lambda = i \mathbb{R}$ of the explicit RK methods and the preconditioned RK methods with TASE operators up to fourth-order. The parameter $\alpha$ in TASE operators is set as the minimum value defined in Table. 1.} \label{Fig_RK_numerical_proof}
\vskip 0.1in
\end{figure}

(b) Analytically investigating whether the stability function defined by Eq.~(\ref{injected_multinomial}) satisfies the condition \RevA{$| \sigma(z) | \leq 1$} for purely imaginary $z$ is difficult. We have shown previously in Section~\ref{sec_alpha} that a necessary condition is given by Eq.~(\ref{eq23_1}), which also applies to the case $z=i y$, with real $y$. However, theorem $\ref{asymptotic_TASE}$ does not guarantee the stability of the method on the entire imaginary axis for intermediate values of the time-step $\Delta t$. Large values of the time step $\Delta t$ are supposedly the worst case for stability, yet nothing guarantees that this holds mathematically. In order to prove that for any $\Delta t$, the time integration of Eq.~(\ref{eq21_2}) remains stable as long as Eq.~(\ref{eq23_1}) is satisfied, and complete the present analysis, we use a computer program to plot $\sigma(i y)$ using a dense sampling of the imaginary axis over very large values of $y$. The results are reported in Fig.~\ref{Fig_RK_numerical_proof}. The plots show consistent results with the analysis of the stability diagrams made in Section~\ref{sec_stability_digram}. In cases where the method is only nearly unconditionally stable, it is worth noting that the detrimental maximal value of the stability function is at most approximately $1.02$.

Although no analytic proof is provided for the second part of corollary \ref{Astability_proof}, the numerical proof performed in Fig.~\ref{Fig_RK_numerical_proof} can be repeated for another time-marching scheme as needed.

\section{\RevB{Extension to general problems and numerical verifications}}
\label{Numerical_validations}

In this section, a set of stiff physical problems, including those with over-resolved spatial resolution, inhomogeneous boundary conditions,  multiple stiffphysical process, nonlinear operators, and grid-induced stiffness, will be conducted to assess the stability and accuracy of the proposed TASE framework.

\subsection{Prototype model examples}  \label{Pedagogical_examples}

Consider the 1D diffusion of a scalar $y(x,t)$ in the infinite domain $x \in \mathbb{R}$ described by the PDE
\begin{equation}
\frac{\partial y}{\partial t} = \frac{\partial ^2 y}{\partial x^2} + A \sin{\left(t/\tau_s\right)}, \label{sec3:eq1_1}
\end{equation}
where $A$ denotes the amplitude of the unsteady source term and $\tau_s=50$ is its time scale. Assuming the initial condition
\begin{equation}
y(x,0)=1 - \cos{(x)}, \label{sec3:eq1_2}
\end{equation}
the exact solution is 
\begin{equation}
y(x,t)=1 - \cos{(x)} \exp{\left(- t\right)} - A \tau_s \left( \cos{\left(t / \tau_s \right)} - 1\right). \label{sec3:eq1_3}
\end{equation}

Since the solution is $2\pi$-periodic, the numerical solution is sought on a spatial grid $0 \leq x \leq 2\pi$ and periodic boundary conditions are used. The grid resolution $\Delta x=2\pi/N$, where $N$ is the total number of collocation points used to describe the solution, should be chosen small enough to resolve the smallest physical length scale of the problem. In this problem, it is the lengthscale of the initial condition, and is approximated as its wavenumber $\delta_{\mathrm{physics}}=1$. The temporal resolution $\Delta t$ must be smaller than the physical time scales of the problem. There are two physical time scales in this problem, one that describes the fast diffusion of the initial spatial inhomogeneities, $\tau_{\mathrm{physics}}=\delta_{\mathrm{physics}}^2/1=1$, and another that represents the slow uniform oscillation imposed by the source, $\tau_s=50$. Depending on the grid resolution $\Delta x$, the final time of interest $\tau_f$, and the source amplitude $A$, an alternative to traditional explicit schemes becomes motivated for one of the reasons described in Section~\ref{sec1}. The objective of this section is to demonstrate the accuracy and stability properties of TASE operators in all of these scenarios, and remind how the choice of the time step remains relevant in these cases. Over-resolved, quasi-steady and steady cases are examined in Sections~\ref{sec31}, \ref{sec32} and \ref{sec33}, respectively. The semi-discrete version of Eq.~(\ref{sec3:eq1_1}) resembles that given in Eq.~(\ref{eq21_1}), where $L$ is the numerical approximation of the second-order derivative operator, but is supplemented with a time dependent source term on the right hand side that is independent from the diffusion term. This second term is not modified in the following, and further discussions and justifications about general treatment of source terms is provided in Section~\ref{boundary_problem}.

In all cases, the maximum time step for numerical stability is given by $\Delta t_{\mathrm{stability}}=C\Delta x^2/D$, where $C$ given in Table~\ref{tab_alpha} depends on the chosen explicit time-marching scheme, and $D$ depends on the numerical approximation of the second-order derivative. For second-order central finite-difference, fourth-order central finite-difference, and Fourier differentiation, $D$ is respectively equal to $4$, $16/3$ and $\pi^2$ \cite{moin2010fundamentals}.

\subsubsection{Scenario 1: the physical features are spatially over-resolved} \label{sec31}

\begin{figure}[t!] 
    \vskip 0.1in
    \centering
    \includegraphics[height=\textwidth,angle=-90,clip]{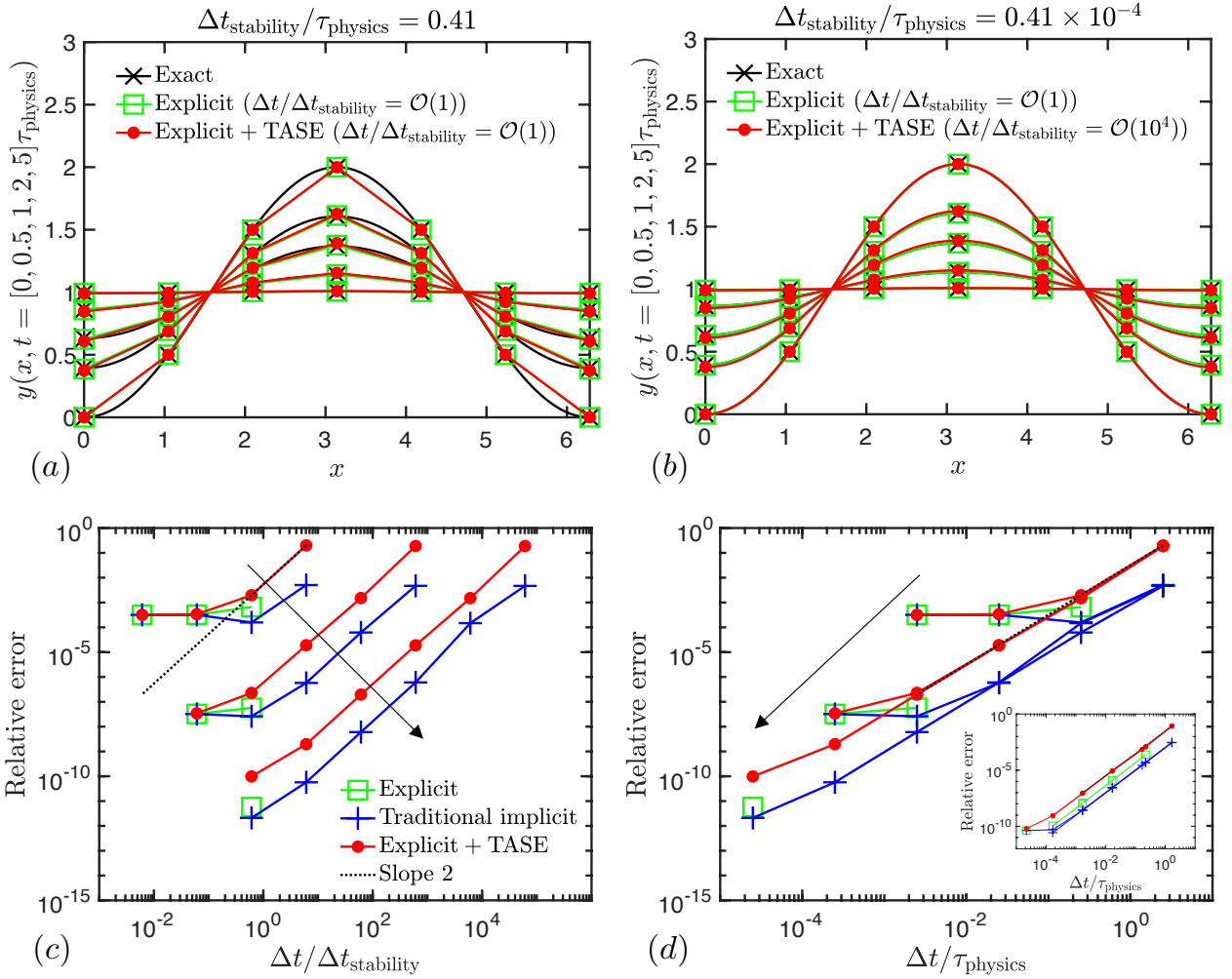}
    \caption{Solution profiles at five different times for (a) the least over-resolved and (b) the most over-resolved case. Relative error plots plotted versus (c) $\Delta t/\Delta t_{\mathrm{stability}}$ and (d) $\Delta t/\tau_{\mathrm{physics}}$. The inset in (d) shows the same figure with Fourier differentiation instead of fourth-order central finite difference. The black arrow in the bottom panel indicates increase over-resolution. \RevB{The compared numerical methods include Explicit (RK2), Traditional implicit (SDIRK2), and Explicit + TASE (RK2 + TASE2).}}\label{fig3}
    \vskip 0.1in
\end{figure}

This section treats the example of over-resolving the spatial scales of the problem, in a way that results in an unphysical constraint on the maximum allowed time step for numerical stability. The short-time solution of Eq.~(\ref{sec3:eq1_1}) at final time is $\tau_f/\tau_{\mathrm{physics}}=5$ is of interest, the source amplitude $A$ is set to 0 and degree of over-resolution is varied with the number of points from 6 to 60 to 600. The spatial and time derivatives are approximated using fourth-order central finite-difference and second-order RK \RevA{(see Table~\ref{But_ERK2} in \ref{But_App})}, respectively. \RevA{The baseline implicit method is a second-order SDIRK method \RevA{(see Table~\ref{But_SDIRK2} in \ref{But_App})}.}

Figure~\ref{fig3}(a,b) show the solution profile at five times during the transient decay of the solution for two cases: the least and the most over-resolved case that respectively correspond to resolving the initial profile with $N=6$ or $600$ points, respectively, which correspond to $\Delta x/\delta_{\mathrm{physics}}=1.05 \times 10^{0}=\mathcal{O}(1)$ and $1.05 \times 10^{-2} = \mathcal{O}(10^{-2})$. In both cases, the time step for the explicit method is equal to the maximum time step for linear stability, with the exception of being rounded so that an integer number of time steps reaches exactly the final time of interest, \textit{i.e.} $\Delta t/\Delta t_{\mathrm{stability}}=6.08 \times 10^{-1}=\mathcal{O}(1)$, which in turn yields $\Delta t/\tau_{\mathrm{physics}}=2.5 \times 10^{-1}=\mathcal{O}(10^{-1})$ and $2.5 \times 10^{-5}=\mathcal{O}(10^{-4})$ for $N=6$ and $600$, respectively. On the contrary, the TASE solution is obtained with a physics-based time step to resolve the smallest physical time scale of the problem, i.e.  $\Delta t/\tau_{\mathrm{physics}}=2.5 \times 10^{-1}=\mathcal{O}(10^{-1})$, irrespectively of the maximum time step of stability. This yield $\Delta t/\Delta t_{\mathrm{stability}}=6.08 \times 10^{-1}=\mathcal{O}(1)$ and $6.08 \times 10^{3}=\mathcal{O}(10^4)$ for $N=6$ and $600$, respectively. Profiles obtained with the different numerical methods are almost indistinguishable to plotting accuracy. The remarkable agreement is obtained while in the second case the explicit + TASE solver uses a time step larger than that used by the explicit solver by four orders of magnitude.

Figure~\ref{fig3}(c) shows the relative error of the final solution profile in function of the time step normalized by the maximum time step allowed for linear stability, $\Delta t/\Delta t_{\mathrm{stability}}$, where the different curves correspond to three different levels of over-resolution $\Delta x/\delta_{\mathrm{physics}}=\mathcal{O}(1)$, $\mathcal{O}(10^{-1})$ and $\mathcal{O}(10^{-2})$. It confirms the expected order of accuracy of the TASE method. With the addition of the TASE operator, the explicit method has preserved its order of accuracy and becomes unconditionally stable. The error plots shifts downward with increasing over-resolution, which shows that the more over-resolved the initial profile is, the more justified the use of a time step larger than the largest one allowed by stability analysis. The error of the TASE method is slightly superior to that of the explicit method without TASE due to the error introduced by the TASE operator, although it remains a small error. The TASE error is larger than that of the baseline implicit method used in this example. The error saturates for the least two over-resolved cases as it becomes dominated by the spatial differentiation error. For each numerical method, all curves for a given time-marching scheme collapse when plotted versus $\Delta t/\tau_{\mathrm{physics}}$ instead, owing to the fact that the only physics-based requirement to obtaining an accurate solution is the resolution of the physical time scale of the initial solution profile. The plot shows that the error approaches order unity when this ratio becomes of order unity. In order to attain a given error within the bounds of engineering interest, the required time step that can be used with the TASE method when the solution is over-resolved is unusable with traditional explicit methods, that would be unstable in that case. In the inset, the same plot is shown when the Fourier differentiation is used instead of fourth-order central finite difference. Since the spatial differentiation becomes exact in this case, the plot reflect time errors only and the the saturation at small time steps disappears.

\begin{table}[!ht]
\begin{center}
\begin{tabular}{ c | c | c | c | c | c | c}
& \multirow{2}{*}{\textbf{RK}} & \multirow{2}{*}{\textbf{SDIRK}} & \multicolumn{4}{| c}{\textbf{RK + TASE}}\\
\cline{4-7}
& & & $\mathbf{p=1}$ & $\mathbf{p=2}$ & $\mathbf{p=3}$ & $\mathbf{p=4}$ \\
\hline
\multirow{2}{*}{\textbf{k=1}} & $5$ & $8$    & $4$ & - & - & - \\
                              & -   & $(33)$ & $(22)$ & - & - & - \\
                            \hline
\multirow{2}{*}{\textbf{k=2}} & $9$ & $17$    & $8$ & $8$ & - & - \\
                              & -   & $(573)$ & $(423)$ & $(959)$ & - & - \\
                            \hline
\multirow{2}{*}{\textbf{k=3}} & $14$ & $31$    & $16$    & $16$     & $15$     & - \\
                              & -    & $(950)$ & $(626)$ & $(1356)$ & $(2376)$ & - \\
                            \hline
\multirow{2}{*}{\textbf{k=4}} & $20$ & $33$    & $20$    & $19$     & $16$     & $16$ \\
                              & -    & $(949)$ & $(941)$ & $(1885)$ & $(3138)$ & $(4165)$
\end{tabular}
\end{center}
\caption{\RevA{CPU time ($10^{-5}s$) per time step for the three methods compared in Sec.~\ref{sec31}, calculated for $N=600$ and a fourth-order spatial differencing scheme. The values in parentheses correspond to the more general implementation that relies on computing the required matrix inversions on the fly, as opposed to the more efficient approach that consists in storing the pre-computed inverse matrices. All SDIRK methods have the same number of stages as their order, except SDIRK4  that has three stages. In cases where the matrix inverses are precomputed, the inverses are computed using the Matlab built-in \textit{inv} function. In cases where the linear systems of equations are solved on the fly, the solutions are obtained using the Matlab built-in \textit{backslash} operator.}} \label{time_table1}
\end{table}

\RevA{Table~\ref{time_table1} shows a comparison of the elapsed time per computational time step for different time integration schemes from second-order to fourth-order accuracy. Considering the linear nature of the operator $L$, the relevant matrix operators can either be precomputed once or on the fly multiple times in practice. With the former, it is shown that the efficiency of present RK+TASE scheme is in general slightly better than that of classical implicit SDIRK scheme with the same accuracy order. Otherwise, the RK+TASE scheme is more expensive than SDIRK, since SDIRK solely needs one matrix inversion per substep whereas the TASE operator involves as many inversions as the accuracy order.}

In this example, over-resolution of the solution justifies the use of a large time step, that cannot be used with an explicit method, as long as it resolves the smallest physical time scale of interest. The more over-resolution, the more useful TASE becomes. In the limit of a normally resolved feature, TASE remains accurate. This is satisfied by construction of the TASE operator, and is a useful characteristic for multi-scale problems that are inhomogeneous and intermittent.

\subsubsection{Scenario 2: the user is interested in the steady-state solution of a stiff problem} \label{sec32}

\begin{figure}[t!] 
    \vskip 0.1in
    \centering    \includegraphics[height=\textwidth,angle=-90,clip]{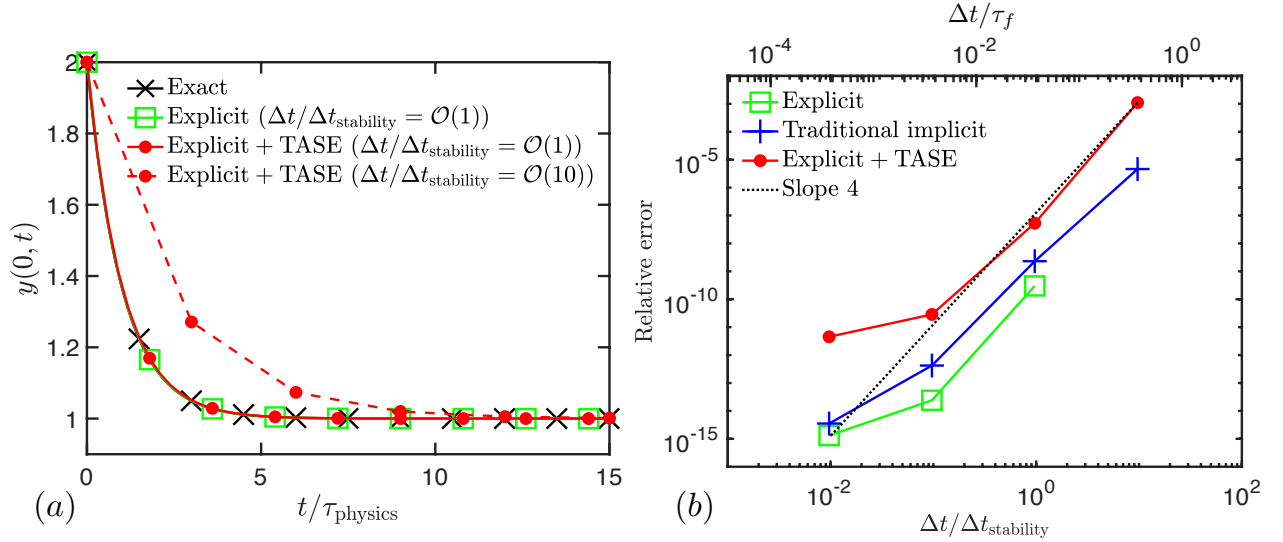} 
    \caption{(a) Time decay of the solution at $x=0$. (b) Relative error plots plotted versus $\Delta t/\Delta t_{\mathrm{stability}}$ (bottom axis) and $\Delta t/\tau_f$ (top axis). \RevB{The compared numerical methods include Explicit (RK4), Traditional implicit (SDIRK4), and Explicit + TASE (RK4 + TASE4).}}\label{fig4}
    \vskip 0.1in
\end{figure}

Now consider the case where $\Delta x / \delta_{\mathrm{physics}} = \mathcal{O}(1)$, i.e. $N=6$, but contrary to before, the user is interested in the final steady-state solution, i.e. $\tau_f \gg \tau_{\mathrm{physics}}$. Here we choose $\tau_f/\tau_{\mathrm{physics}}=15$. The spatial and time derivatives are approximated using exact Fourier spectral differentiation and fourth-order RK \RevA{(see Table~\ref{But_ERK4} in \ref{But_App})}, respectively. \RevA{The baseline implicit method is a fourth-order SDIRK method (see Table~\ref{But_SDIRK4} in \ref{But_App}).}

Figure~\ref{fig4}(a) show the time decay of the solution at $x=0$. The time step for the explicit method is approximately equal to the maximum time step for linear stability, i.e. $\Delta t/\Delta t_{\mathrm{stability}}=8.07 \times 10^{-1} = \mathcal{O}(1)$, which corresponds to $\Delta t/\tau_f=1.67 \times 10^{-2} = \mathcal{O}(10^{-2})$. On the contrary, the TASE methods need not resolve this time scale but only the time scale of interest, which is here the final time $\tau_f$. For illustration, Figure~\ref{fig4}(a) shows the profile obtained with the TASE methods for $\Delta t/\Delta t_{\mathrm{stability}}=8.07 \times 10^{-1} = \mathcal{O}(1)$ and $8.07 \times 10^{0} = \mathcal{O}(10)$, which correspond to $\Delta t/\tau_f=1.67 \times 10^{-2} = \mathcal{O}(10^{-2})$ and $1.67 \times 10^{-1} = \mathcal{O}(10^{-1})$, respectively. Profiles obtained with the different numerical methods are almost indistinguishable to plotting accuracy in the limit of small time step. However, the prediction of the solution at the final time of interest is predicted remarkably well by all schemes irrespective of the chosen time step, because in all cases the physics-based constraint $\Delta t/\tau_f \ll 1$ is satisfied.

Figure~\ref{fig4}(b) shows the relative error of the final solution profile in function of the time step normalized by the maximum time step allowed for linear stability, $\Delta t/\Delta t_{\mathrm{stability}}$, and normalized by the final time of interest $\Delta t/\tau_f$. It confirms the expected order of accuracy of the TASE method. Compared to the previous example, it shows the flexibility of TASE to handle different order. With the addition of the TASE operator, the explicit method has preserved its order of accuracy and become unconditionally stable. The plot shows that the error approaches order unity when the time step becomes of the order of the final time step of interest. In order to predict a quantity of interest within error bounds of engineering interest, the required time step that can be used with the TASE method is unusable with traditional explicit methods, that would be unstable in that case.

In this example, the separation of scale between the smallest physical time scale that corresponds to the initial decay of the solution and the final time scale of interest justifies the use of a time step much larger than the maximum time step predicted by linear stability analysis. The larger the separation of scale, the more useful TASE becomes. Here, the separation was by a factor fo 10 only, but in problems of practical interest, it can easily be much larger. 

\subsubsection{Scenario 3: the stiff physics is quasi-steady and irrelevant to the long term transient solution} \label{sec33}

\begin{figure}[t!] 
    \vskip 0.1in
    \centering
    \includegraphics[height=\textwidth,angle=-90,clip]{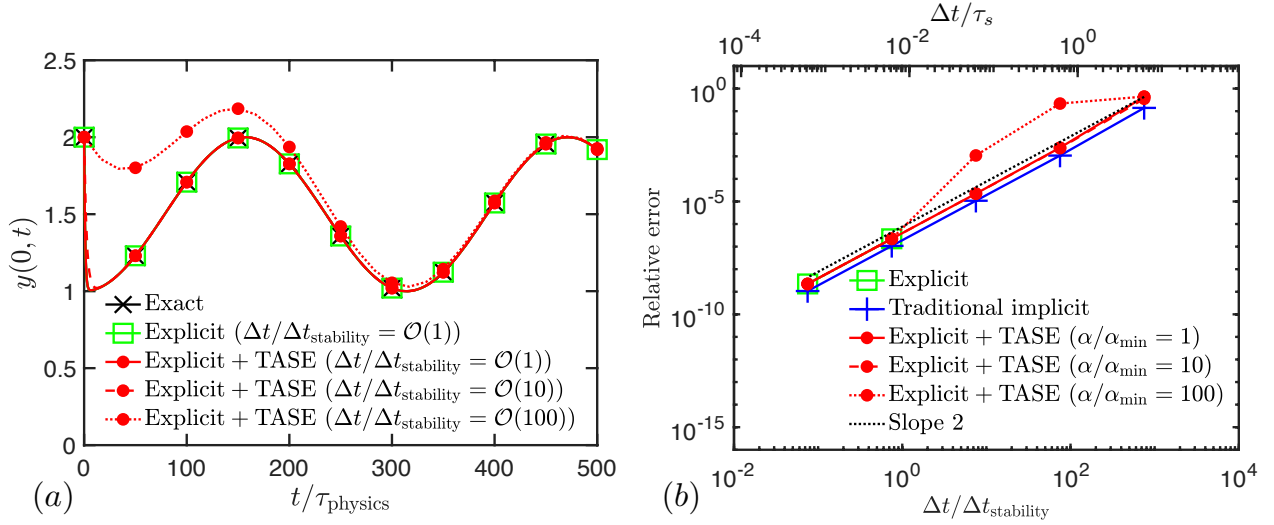} 
    \caption{(a) Time decay of the solution at $x=0$. (b) Relative error plots plotted versus $\Delta t/\Delta t_{\mathrm{stability}}$ (bottom axis) and $\Delta t/\tau_s$ (top axis). \RevB{The compared numerical methods include Explicit (RK2), Traditional implicit (SDIRK2), and Explicit + TASE (RK2 + TASE2).} }\label{fig5}
    \vskip 0.1in
\end{figure}

Now consider the case where $\Delta x / \delta_{\mathrm{physics}} = \mathcal{O}(1)$, as in Section~\ref{sec32} but the source term is activated, and the user is not interested in the initial quasi-steady fast decay of the initial profile. Instead the objective is only to accurately capture the slow oscillations encoded by the source term, i.e. $\tau_\mathrm{physics} \ll \tau_\mathrm{physics}^{NQS} =\tau_s$. The amplitude source is $A=0.01$. The spatial and time derivatives are approximated using exact Fourier spectral differentiation and second-order RK, respectively. \RevA{The baseline implicit method is a second-order SDIRK method \RevA{(see Table~\ref{But_SDIRK2} in \ref{But_App})}.}

Figure~\ref{fig5}(a) show the time development of the solution at $x=0$. The time step for the explicit method is equal to the maximum time step for linear stability, i.e. $\Delta t/\Delta t_{\mathrm{stability}}=7.5 \times 10^{-1}=\mathcal{O}(1)$, which in turn yields $\Delta t/\tau_s=3.3 \times 10^{-3}=\mathcal{O}(10^{-3})$. On the contrary, the TASE methods need not resolve this time scale but only the time scale of interest, which is here the time scale of the source $\tau_s$. For illustration, Figure~\ref{fig5}(a) shows the profile obtained with the TASE methods for $\Delta t/\Delta t_{\mathrm{stability}}=7.5 \times 10^{-1}=\mathcal{O}(1)$, $7.5=\mathcal{O}(10)$ and $75=\mathcal{O}(10^2)$, which correspond to $\Delta t/\tau_s=3.3 \times 10^{-3}=\mathcal{O}(10^{-3})$, $3.3 \times 10^{-2}=\mathcal{O}(10^{-2})$ and $3.3 \times 10^{-1}=\mathcal{O}(10^{-1})$, respectively. Profiles obtained with the different numerical methods are almost indistinguishable to plotting accuracy in the limit of small time step. However, the prediction of the solution after an initial fast transient decay is predicted remarkably well by all schemes irrespective of the chosen time step, because in all cases the physics-based constraint $\Delta t/\tau_s \ll 1$ is satisfied.

Figure~\ref{fig5}(b) shows the relative error of the final solution profile in function of the time step normalized by the maximum time step allowed for linear stability, $\Delta t/\Delta t_{\mathrm{stability}}$, and normalized by the source time scale $\Delta t/\tau_s$. It confirms the expected order of accuracy of the TASE method. With the addition of the TASE operator, the explicit method has preserved its order of accuracy and become unconditionally stable. The plot shows that the error approaches order unity when the time step becomes of the order of the source time scale. In order to predict a given quantity of interest within the error bounds of engineering interest, the required time step that can be used with the TASE method is unusable with traditional explicit methods, that would be unstable in that case. The effect of $\alpha$ is also shown. With increasing $\alpha$, the order of convergence is preserved, but the error is increased. Note that in practice it is not a tunable parameter, as theoretically prescribed as described in Section~\ref{TASE_operators}.

In this example, the separation of scale between the smallest physical time scale that corresponds to the initial decay of the solution and the source time scale justifies the use of a time step much larger than the maximum time step predicted by linear stability analysis. The larger the separation of scale, the more useful TASE becomes. Here, the separation was by a factor of 100 only, but in problems of practical interest, it can easily be much larger.

\subsection{Application to problems with inhomogeneous boundary conditions (source terms)} \label{boundary_problem}

\subsubsection{Method} \label{sec41}

In this section, we treat a more general example where the boundary conditions are implemented as a source term in the semi-discrete equation. In this case, the semi-discrete system typically takes the form
\begin{equation}
\frac{d Y}{d t} = L Y + S, \label{eq41_1}
\end{equation} 
where $S$ is a vector that encodes for the physical boundary conditions, and that may depend on $L$ and the specific choice of numerical operators used at the boundaries. When the TASE methodology is used, the TASE operator simply needs to be applied to both terms simultaneously, yielding the modified semi-discrete equation
\begin{equation}
\frac{d Y}{d t} = T^{(p)}_L \left( L Y + S \right), \label{eq41_2}
\end{equation}
contrary to the inaccurate one that would be obtained if the TASE operator was only applied to the linear operator, $dY/dt = T^{(p)}_L L Y + S$. The simultaneous application of the TASE operator to both terms is necessary since both encode for a single physical process and compete near the boundary. In particular, this approach guarantees that the correct quasi-steady or steady solution of the system is captured, as can be understood by substituting $dY/dt=0$ in Eq.~(\ref{eq41_2}), and solving for $Y$. The presence of a source term originating from the boundary condition discretization is a particular case of the more general case of competing operators, which is treated in Section~\ref{sec51}.

\subsubsection{Example} \label{sec42}

\begin{figure}[t!] 
    \vskip 0.1in
    \centering
    \includegraphics[height=\textwidth,angle=-90,clip]{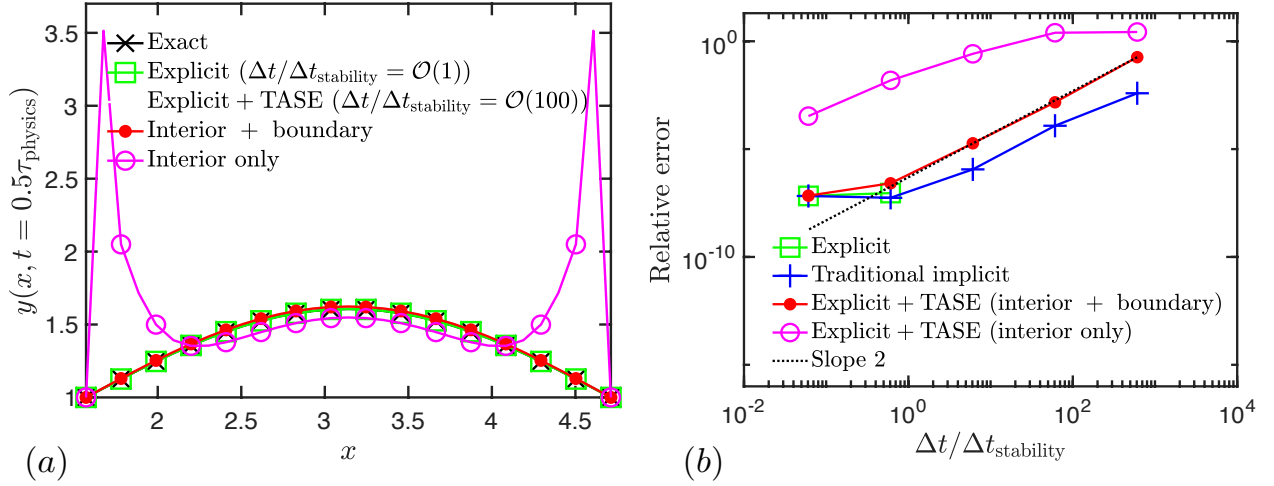} 
    \caption{(a) Solution profile at $t=0.5 \tau_{\mathrm{physics}}$. (b) Relative error plots plotted versus $\Delta t/\Delta t_{\mathrm{stability}}$. \RevB{The compared numerical methods include Explicit (RK2), Traditional implicit (Crank-Nicolson), and Explicit + TASE (RK2 + TASE2).}}\label{fig6}
    \vskip 0.1in
\end{figure}

In order to demonstrate the method for handling source term, and how the failure manifests if the boundary conditions implemented as a source term is not subject to the TASE operator, the example from Section~\ref{sec31} is repeated in the subdomain $\pi/2 \leq x \leq 3\pi/2$ so that the exact boundary Dirichlet boundary conditions $y(\pi/2,t)=y(3\pi/2,t)=1$ yields a non-zero source term in the semi-discrete equation. The example uses $N=30$ collocation points, and the spatial and time derivatives are approximated using fourth-order central differencing and second-order RK schemes \RevA{(see Table~\ref{But_ERK2} in \ref{But_App}), respectively. The baseline implicit method is the Crank-Nicolson method (see Table~\ref{But_CN} in \ref{But_App}). The Dirichlet boundary conditions are locally enforced using a second-order finite difference scheme.}

Figure~\ref{fig6}(a) shows a comparison of the instantaneous solution profile obtained at an intermediate time $t=0.5\tau_{\mathrm{physics}}$ for different numerical methodologies. The time step for the explicit method is equal to the maximum time step for linear stability, i.e. $\Delta t/\Delta t_{\mathrm{stability}}=6.08 \times 10^{-1} = \mathcal{O}(1)$. On the contrary, the TASE methods use a physics-based time step to resolve the smallest physical time scale of the problem, i.e.  $\Delta t/\tau_{\mathrm{physics}}=2.50 \times 10^{-1} = \mathcal{O}(10^{-1})$, irrespectively of the maximum time step of stability. This yield $\Delta t/\Delta t_{\mathrm{stability}}=6.08 \times 10^{1} = \mathcal{O}(10^2)$. The figures compares the profiles obtained with the correct implementation of the TASE methodology, where both linear operator and source term are modified, to the incorrect one where only the linear operator is modified, as described in Sec.~\ref{sec41}, with $\Delta t/\Delta t_{\mathrm{stability}}=100$. The figure shows that similar conclusions as in Sec.~\ref{sec31} applies when the TASE methodology is correctly implemented. On the contrary, large errors arise near the boundaries when the boundary source term is not modified, and subsequently propagate in the interior of the domain.
 
The observations described above are confirmed in Fig.~\ref{fig6}(b) that shows the relative error of the final solution profile in function of the time step normalized by the maximum time step allowed for linear stability, $\Delta t/\Delta t_{\mathrm{stability}}$. Besides conclusions already drawn in previous part of the paper, it is worth noting that the application of TASE operators with wrong handling of boundary source terms preserves its order of accuracy but much larger errors due to the inconsistent numerical treatment of the diffusion physics in the interior and at the edges of the domain.

\subsection{Application to problems with multiple stiff physical processes (operator splitting)} \label{stiff_problem}

\subsubsection{Method} \label{sec51}

In this section, we treat a more general example where the problem involves multiple stiff linear operators $L_1$ and $L_2$
\begin{equation}
\frac{d Y}{d t} = L_1 Y + L_2 Y. \label{eq51_1}
\end{equation}

The TASE operators can be used to bypass the linear stability restrictions of each operators. Each operator can be modified combined
\begin{equation}
\frac{d Y}{d t} = T^{(p)}_{L_1+L_2} \left( L_1 Y + L_2 Y \right), \label{eq51_2}
\end{equation}

or separately

\begin{equation}
\frac{d Y}{d t} = T^{(p)}_{L_1} L_1 Y + T^{(p)}_{L_2} L_2 Y. \label{eq51_3}
\end{equation} 

The choice of one strategy versus the other is problem-dependent, and depends on whether the underlying quasi-steady process to be bypassed by the temporal resolution involves competition of both operators. When the operators compete, i.e. they are both simultaneously active and significant in a (quasi-)steady process at the same location, a combined TASE operator $T^{(p)}_{L_1 + L_2}$ must be constructed and applied to both terms, as in Eq.~(\ref{eq51_2}). Alternatively, when they do not compete, it is sufficient to build distinct TASE operators $T^{(p)}_{L_1}$ and $T^{(p)}_{L_2}$ and apply them to the respective terms in the equation, as in Eq.~(\ref{eq51_3}). This physics-based constraint for splitting or grouping the application of TASE to the right-hand side operators is reminiscent of the constraint that underpins \RevA{the use of operator splitting methods} \cite{marchuk1968some,strang1968construction}. In general, competing operators should not be split but computed simultaneously in this latter methodology \cite{strang1968construction}. In the absence of knowledge about the competition among physical operators, the user may safely adopt the strategy described in Eq.~(\ref{eq51_2}). It is more general but since its implementation can be more cumbersome and more costly than the alternative one, it is desirable to use Eq.~(\ref{eq51_3}) when prior knowledge about the non-competing stiff operators exists.

\subsubsection{Example} \label{sec52}

\begin{figure}[t!] 
    \vskip 0.1in
    \centering
    \includegraphics[height=\textwidth,angle=-90,clip]{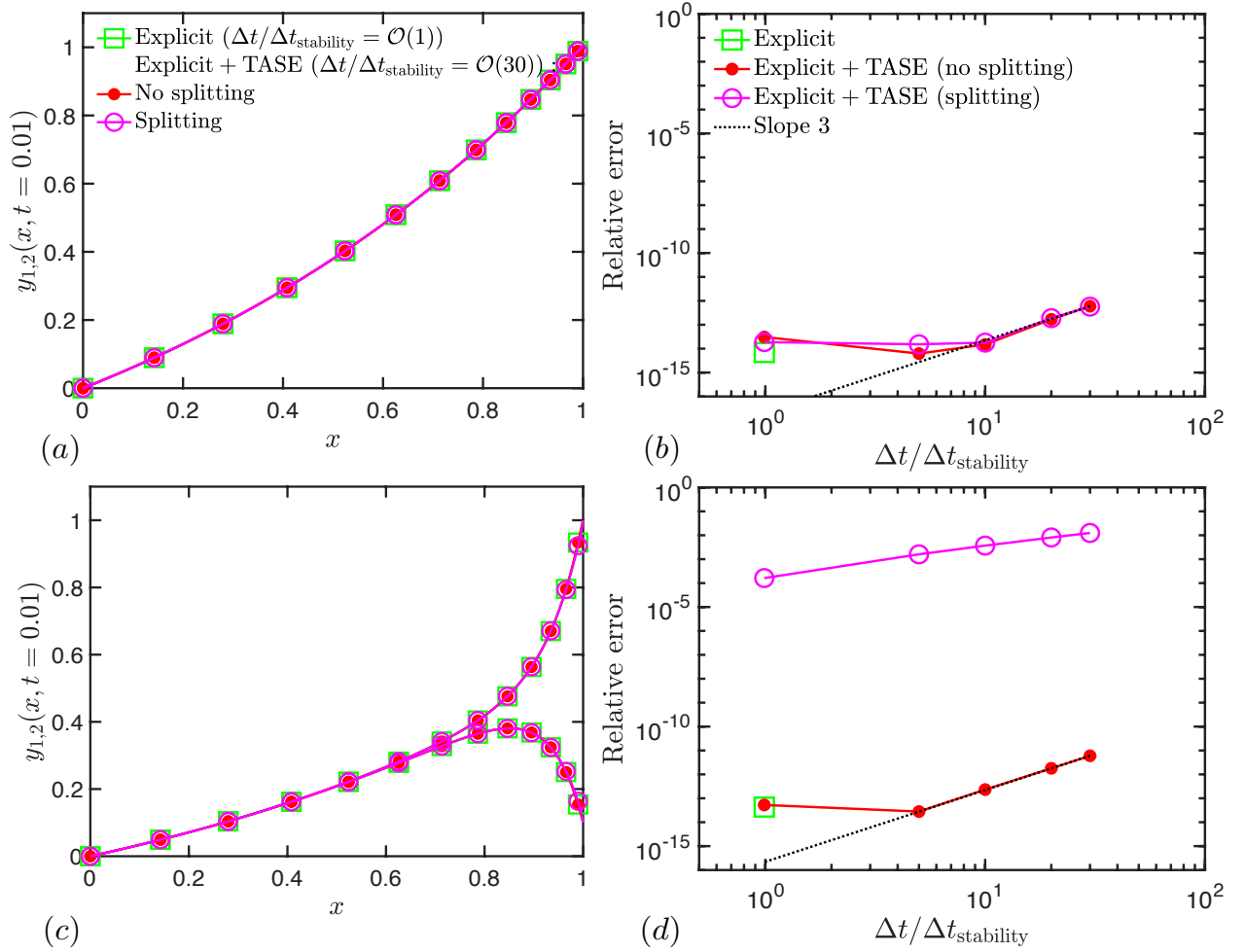}
    \caption{(a,c) Solution profile at $t=10^{-2}$. (b,d) Relative error plots plotted versus $\Delta t/\Delta t_{\mathrm{stability}}$. The top panel corresponds to $y_1(x=1,t)=y_2(x=1,t)=1$, while the bottom panel corresponds to $y_1(x=1,t)=1$ and $y_2(x=1,t)=0.1$. \RevUs{The relative error is computed using the infinite norm to emphasize the local origin of the error near $x=1$ in the bottom panel for the splitting case.} \RevB{The compared numerical methods include Explicit (RK3) and Explicit + TASE (RK3 + TASE3).}}\label{fig7}
    \vskip 0.1in
\end{figure}

Consider the two-species advection-diffusion-reaction
\begin{equation}
\frac{\partial y_1}{\partial t} + U \frac{\partial y_1}{\partial x} = D \frac{\partial^2 y_1}{\partial x^2} + K(y_2-y_1), \label{eq}
\end{equation}
\begin{equation}
\frac{\partial y_2}{\partial t} + U \frac{\partial y_2}{\partial x} = D \frac{\partial^2 y_2}{\partial x^2} - K(y_2-y_1), \label{eq52_1}
\end{equation} 
in the bounded domain $0 \leq x \leq 1$. The semi-discrete equation becomes
\begin{equation}
\frac{dY}{d t} = \left( L_c + L_d +L_r \right) Y + \left( S_c+S_d \right), \label{eq52_2}
\end{equation} 
where $Y=[Y_1 Y_2]^{T}$ is the concatenated vector of unknowns, $L_c$, $L_d$, and $L_r$ denote the discretized operators that correspond to the convection, diffusion, and reaction terms, respectively. The source terms are denoted $S_c$ and $S_d$, respectively. The numerical algorithm employs a non-uniform mesh stretched near the right-boundary in order to capture the fine feature that develops over time there, and a finite-volume discretization of the equation that employs second-order spatial derivative operators, with $N=50$ control volumes. The time-marching scheme is \RevUs{a third-order Runge-Kutta method} \RevA{(see Table~\ref{But_ERK3} in \ref{But_App})}. The problem parameters $U$, $D$ and $K$ are chosen in such a way that the reaction time scale is much faster that the convection and diffusion one, i.e. $K=10^4$, $D=10^2$ and $U=10^2$. Assuming that the physical length scale is $\delta_{\mathrm{physics}}=D/U=1$, this renders the following correspond time scales: $\tau_r=1/K=10^{-4}$, $\tau_c = \delta_{\mathrm{physics}} / U = 10^{-2}$, and $\tau_d = \delta_{\mathrm{physics}}^2 / D = 10^{-2}$, i.e $\tau_r \ll \tau_{c}, \tau_{d}$. The restricted time step for each physical operator is given by linear analysis. The final time is chosen to be $\tau_f=10^{-2}$. The maximum time step for linear stability is estimated for each operator separately. It yields $\tau_f/\Delta t=6.6 \times 10^{4}$.

In the following, depending on the choice of boundary conditions, we illustrate when TASE splitting is appropriate. In all cases, the boundary conditions at $x=0$ is fixed for both species $y_1(x=0,t)=y_2(x=0,t)=0$.  On the right at $x=1$, we simulate two different cases. In the first case, both species have the same boundary conditions $y_1(x=1,t)=y_2(x=1,t)=1$, and the initial conditions are $y_1(x,t=0)=x$ and $y_2(x,t=0)=x^2$. In this case, at short time the reaction term is active until $y_1=y_2$ then it vanishes and the steady solution is slowly reached as a result of competition between convection and diffusion. In this case, the appropriate TASE treatment of the semi-discrete equation (\ref{eq52_2}) is one where the reaction term can be split and modified separately from the other convection and diffusion terms, respectively
\begin{equation}
\frac{dY}{d t} = T_{L_c + L_d}^{(p)}\left[\left( L_c + L_d \right) Y + \left( S_c+S_d \right)\right] + T_{L_r}^{(p)}L_r Y, \label{eq52_3}
\end{equation} 
where the source terms are treated as advocated in Sec.~\ref{sec41}. The splitting is achieved by decoupling the physical operators in time domain. In the second case, the two species have different boundary conditions $y_1(x=1,t)=1$ and $y_2(x=1,t)=0.1$, and the initial conditions are $y_1(x,t=0)=x$ and $y_2(x,t=0)=0.1x^2$. The time dynamics is similar as before except that the incompatible boundary condition on the right of the domain render the reactive term active at all times, therefore it permanently competes with convection and diffusion and cannot be modified separately. In this case, the correct treatment is given by 
\begin{equation}
\frac{dY}{d t} = T_{L_c + L_d+L_r}^{(p)}\left[\left( L_c + L_d + L_r \right) Y + \left( S_c+S_d \right)\right]. \label{eq52_4}
\end{equation} 

Figure~\ref{fig7}(a,c) show the instantaneous solution profile obtained with the different approaches in both cases. In the first case, where splitting is physically motivated, the figure confirms that both approaches are valid and splitting yield a solution with very good accuracy properties. This is not the case for the second case, where the solution obtained with splitting differs significantly from the correct solution. Figure~\ref{fig7}(c) shows that the error is largest near the right boundary where the reacting term is non-zero and actively compete with convection-diffusion. Figure~\ref{fig7}(b,d) show the relative error of the final solution profile in function of the time step normalized by the maximum time step allowed for linear stability, $\Delta t/\Delta t_{\mathrm{stability}}$. It confirms the conclusions made above. In all cases, the formal order of accuracy of the TASE solution is not affected but large errors occur when the application of the TASE operator is split but not physically motivated.  Note that the deviation from the theoretical slope is due to the fact that the reference solution used to compute the error is one obtained from using the explicit method at $\Delta t / \Delta t_{\mathrm{stability}}=\mathcal{O}(10^{-1})$, yet the TASE error collapses with the of the explicit method at the smallest time step shown on the figure.

\subsection{Application to non-linear problems} \label{nonlinear_problem}

\subsubsection{Method} \label{sec61}

In the present section, the TASE methodology for non-linear problems is described. The semi-discrete equation under consideration is now
\begin{equation}
\frac{d Y}{d t} = N[Y], \label{eq61_1}
\end{equation} 
where $N$ is a non-linear operator. Application of the TASE methodology to the above equation relies on treating the non-linear term as a source term, and premultiply the right-hand side with the TASE operator $T^{(p)}_{L_{Y}}$ to yield
\begin{equation}
\frac{d Y}{d t} = T_{L_{Y}}^{(p)}N[Y], \label{eq61_2}
\end{equation}
where the linear operator $L_{Y}$ is a linear approximation of the nonlinear operator $N$ in the surrounding of $Y$, formally defined as
\begin{equation}
L_{Y}[(\cdot)] = \left. \left. \left. \frac{\partial}{\partial \epsilon} \right( N\left[Y + \epsilon (\cdot)\right] \right) \right|_{\epsilon=0}. \label{eq61_3}
\end{equation} 

It is worthwhile to note that the above method is equivalent to decomposing the non-linear operator $N$ into its linearized approximation and a remainder, $N[Y] = L_{Y^*}Y + \left( N[Y] - L_{Y^*}Y \right)$, and applying the method described in Sec.~\ref{sec41} to handle linear problems with source terms. The resulting semi-discrete equation collapses to that given in Eq.~(\ref{eq61_2}).

The TASE methodology can be readily applied to non-linear problems and enables to achieve arbitrary high-order accuracy even for non-linear problems, without the need to solve a non-linear system of equations. For a particular problem, the estimation of the linearized approximation of the non-linear operator can be obtained a priori and inputted in the solver, and evaluated at every time-step using the current value of the vector $Y$.

\subsubsection{Example 1} \label{sec62}

\begin{figure}
    \vskip 0.1in
    \centering
    \includegraphics[height=\textwidth,angle=-90,clip]{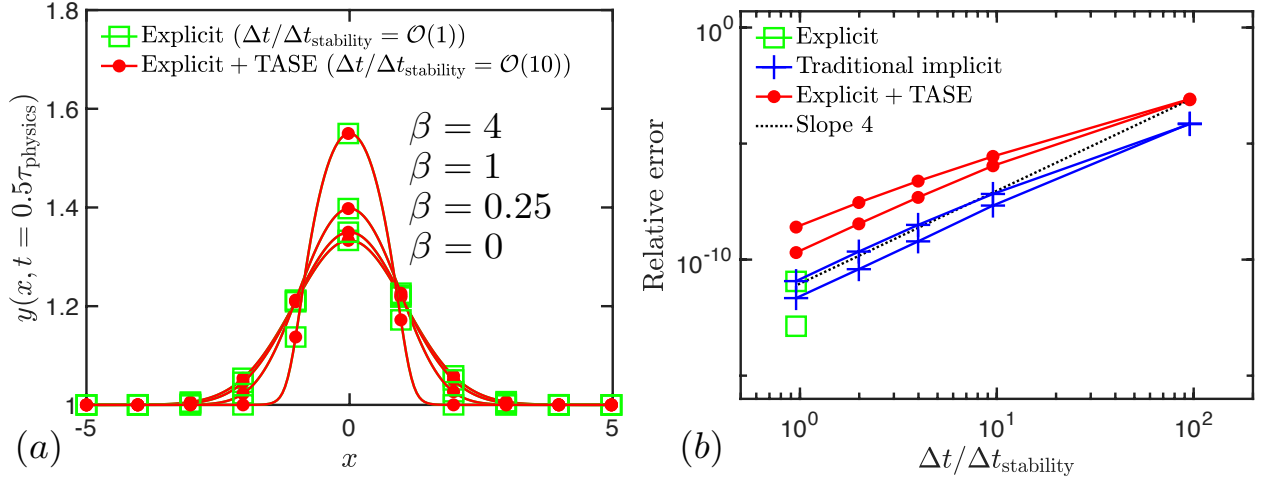} 
    \caption{(a) Solution profiles at $t=0.5\tau_{\mathrm{physics}}$ for different values of $\beta$. (b) Relative error plots plotted versus $\Delta t/\Delta t_{\mathrm{stability}}$ for $\beta=0$ and $\beta=4$. The error curves shift upward with larger values of $\beta$. \RevB{The compared numerical methods include Explicit (RK4), Traditional implicit (SDIRK4), and Explicit + TASE (RK4 + TASE4).}}\label{fig8}
    \vskip 0.1in
\end{figure}

In order to demonstrate the TASE method for handling non-linear semi-discrete equation, consider the PDE
\begin{equation}
\frac{\partial y}{\partial t} = \frac{\partial}{\partial x}\left(\left(\frac{y}{2}\right)^\beta \frac{\partial y}{\partial x}\right) \label{eq62_1}
\end{equation}
where positive non-zero values of $\beta$ renders Eq.~(\ref{eq62_1}) analogous to the 1D diffusion of a scalar with a power-law dependence of the diffusion coefficient. As prescribed by the method described in Sec.~\ref{sec61}, the semi-discrete equation resulting from spatial discretizing Eq.~(\ref{eq62_1}) is premultiplied by the TASE operator $T^{(p)}_{L_Y}$ where $L$ is estimated as 
\begin{equation}
L_{Y}(\cdot) = \frac{\delta}{\delta x}\left(\left(\frac{Y}{2}\right)^\beta \frac{\delta (\cdot)}{\delta x}\right)  + \frac{\delta}{\delta x}\left( \frac{\beta}{2} \left(\frac{Y}{2}\right)^{\beta-1} \frac{\delta Y}{\delta x} (\cdot) \right), \label{eq62_2}
\end{equation}
where $\delta (\cdot)/\delta x$ is used instead of $\partial (\cdot)/\partial x$ to refer to the discrete version of the differentiation operator.

The PDE is numerically solved until $\tau_f=1$ using a finite-volume formulation that uses second-order spatial differencing schemes, and fourth-order RK scheme \RevA{(see Table~\ref{But_ERK4} in \ref{But_App})} in time. \RevA{The results are compared with a fourth-order SDIRK method \RevA{(see Table~\ref{But_SDIRK4} in \ref{But_App})}.} The initial condition is $y(x,t=0)=1+\exp{(-0.25 x^2)}$ and Neumann boundary conditions are used at the edges of the domain $x=\pm 5$. The spatial domain is discretized using $N=200$ control volumes. The errors are computed using the final solution profile obtained with the explicit method using a fine time step.

Figure~\ref{fig8}(a) shows a comparison of the instantaneous solution profile obtained at an intermediate time $t=0.5$ for different numerical methodologies. The time step for the explicit method is equal to the maximum time step for linear stability, i.e. $\Delta t/\Delta t_{\mathrm{stability}}=9.56 \times 10^{-1} = \mathcal{O}(1)$. On the contrary, the TASE methods use a physics-based time step to resolve the smallest physical time scale of the problem, i.e.  $\Delta t/\tau_{\mathrm{physics}}=1.67 \times 10^{-2} = \mathcal{O}(10^{-2})$, irrespectively of the maximum time step of stability, where the physical length scale is estimated using the initial condition as $\delta_{\mathrm{physics}}=1$ and its time scale of decay $\tau_{\mathrm{physics}}=1$. This yield $\Delta t/\Delta t_{\mathrm{stability}}=9.56 \times 10^{0} = \mathcal{O}(10)$. The figure shows that the results obtained with TASE remain stable and accurate to plotting accuracy for a range of exponent $\beta$ ranging from 0 to 0.25 to 1 to 4. Figure~\ref{fig8}(b) confirms that the TASE method preserved the order of accuracy.

%

\begin{table}[!ht]
\begin{center}
\begin{tabular}{ c | c | c | c | c | c | c}
& \multirow{2}{*}{\textbf{RK}} & \multirow{2}{*}{\textbf{SDIRK}} & \multicolumn{4}{| c}{\textbf{RK + TASE}}\\
\cline{4-7}
& & & $\mathbf{p=1}$ & $\mathbf{p=2}$ & $\mathbf{p=3}$ & $\mathbf{p=4}$ \\
\hline
\multirow{1}{*}{\textbf{k=1}} & $0.6$ & $237$    & $2.1$ & - & - & - \\
                            \hline
\multirow{1}{*}{\textbf{k=2}} & $0.9$ & $399$    & $3.3$ & $6.0$ & - & - \\
                            \hline
\multirow{1}{*}{\textbf{k=3}} & $1.2$ & $598$    & $4.6$    & $8.7$     & $11$     & - \\
                            \hline
\multirow{1}{*}{\textbf{k=4}} & $1.5$ & $1303$    & $12$    & $20$     & $30$     & $40$ \\
\end{tabular}
\end{center}
\caption{\RevA{CPU time ($10^{-3}s$) per time step for the three methods compared in Sec.~\ref{sec62}, calculated for $\beta=4$. All SDIRK methods have the same number of stages as their order, except SDIRK4  that has three stages. For the implicit cases, the nonlinear systems of equations are solved using the Matlab built-in \textit{fsolve} function with default parameters.}} \label{time_table2}
\end{table}

\RevA{Table~\ref{time_table2} shows the statistics of the elapsed time per computational time step for different time integration schemes from second-order to fourth-order accuracy. It is shown that the present RK+TASE schemes are more efficient than the classical SDIRK schemes without sacrificing the accuracy order. The rationale behind is that RK+TASE relies on solving a couple of linear systems instead of the expensive nonlinear system required by the classical SDIRK methods. Considering that the RK+TASE method allows for considerably larger time-steps, the overall cost of the present implicit method is also significantly smaller than the standard explicit RK method when differences in absolute time steps are factored in.}
%

\subsubsection{Example 2} \label{sec:ODE}

\begin{figure}
    \vskip 0.1in
    \centering
    \includegraphics[height=\textwidth,angle=-90,clip]{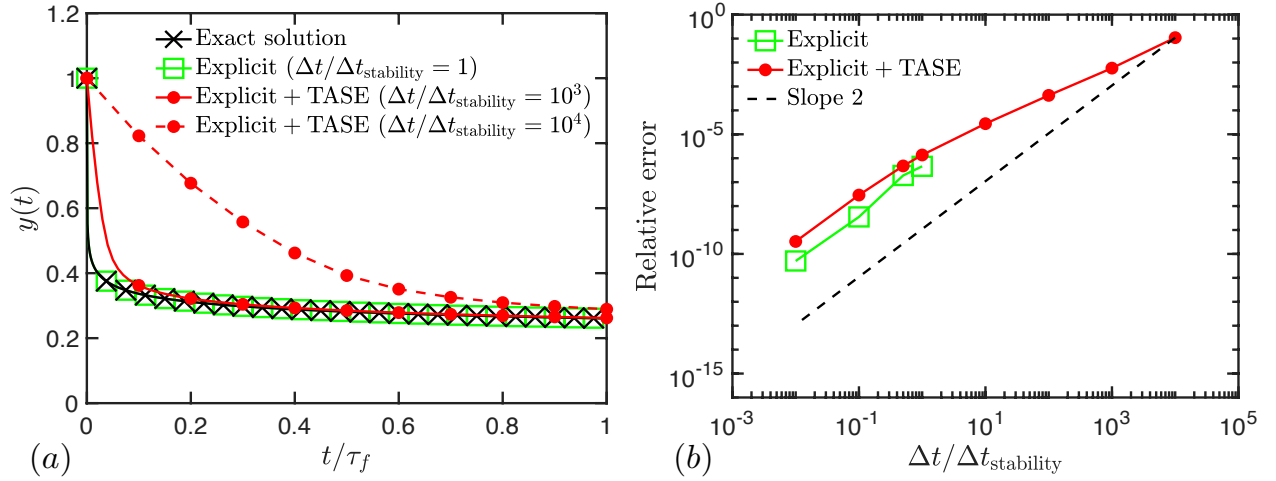} 
    \caption{(a) Time development of the solutions. (b) Relative error plots plotted versus $\Delta t/\Delta t_{\mathrm{stability}}$. \RevB{The compared numerical methods include Explicit (RK2) and Explicit + TASE (RK2 + TASE2).}}\label{fig:ODE}
    \vskip 0.1in
\end{figure}

We further demonstrate the capability of the TASE methodology to handle nonlinear problems with the following nonlinear ODE
\begin{equation}
\frac{dy}{dt} = - y^{\beta}  \label{eq_ODE}
\end{equation}
with $\beta=10$ and initial condition $y(0)=1$. Following the method described in Sec.~\ref{sec61}, the TASE operator $T^{(p)}_{L_y}$ is computed based on the linearized operator $L$ estimated as
\begin{equation}
L_{y}(\cdot) = -\beta y^{\beta - 1}, \label{eq_ODE_L}
\end{equation}
which corresponds to $\beta$ times the operator that would be naively obtained by estimating $N[y]/y$.

The PDE is numerically solved until $\tau_f=2\times 10^4$ using a second-order RK scheme \RevA{(see Table~\ref{But_ERK2} in \ref{But_App})}. The relative error is computed using the solution at the final time and the exact solution $y_{\mathrm{exact}}(t) = \left( 1 + (\beta-1)t \right)^{1/(1-\beta)}$.

Figure~\ref{fig:ODE}(a) shows the time development of the solutions obtained with different numerical methodologies, against the exact solution. The time step for the explicit method is equal to the maximum time step for linear stability, i.e. $\Delta t/\Delta t_{\mathrm{stability}}=1$, where $\Delta t_{\mathrm{stability}} = 2/(\beta) = 0.2$. On the contrary, the TASE method use $\Delta t/\Delta t_{\mathrm{stability}}=10^3$ and $10^4$. The number of time steps to reach the solution at the final time is reduced from $10^5$ to $10^2$ and $10$, respectively. The figure shows that the results obtained with TASE remain stable and accurate to obtain the target long-time solution. Figure~\ref{fig:ODE}(b) confirms that the TASE method preserved the order of accuracy of the original RK scheme.

\subsection{Application to two-dimensional problems in polar coordinates} \label{polar_problem}

\begin{figure}[t!] 
    \vskip 0.1in
    \centering
    \includegraphics[height=\textwidth,angle=-90,clip]{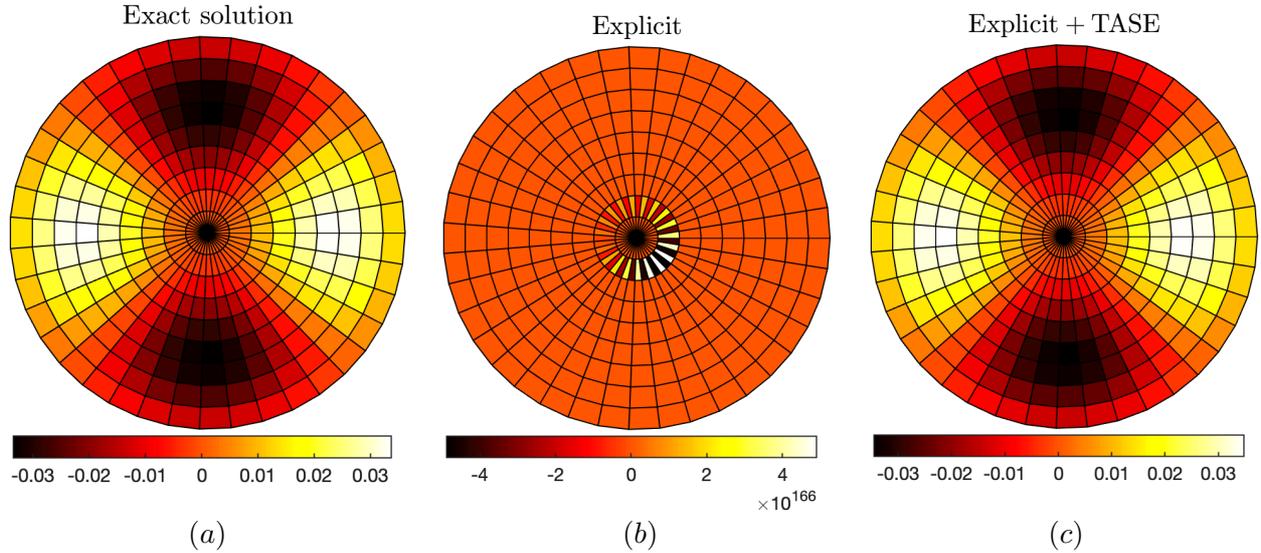} 
    \caption{(a) Exact solution at $t=0.1$, compared with (b) Explicit and (c) Explicit + TASE solution. Both numerical solutions use a time step $\Delta t = 50 \Delta t_{\mathrm{stability,\theta}}$. \RevB{The compared numerical methods include Explicit (RK2) and Explicit + TASE (RK2 + TASE2).}}\label{fig9}
    \vskip 0.1in
\end{figure}

This example demonstrates the applicability of TASE operators to two-dimensional problems formulated in polar coordinates. In this case, the azimuthal over-resolution near the pole, imposed by the need to resolved azimuthal feature far from the pole imposes stringent numerical resolution and stability criterion that can be alleviated by the use of the TASE operator. The numerical stiffness imposed in the near-pole region is a well-known numerical issue \cite{mohseni2000numerical}. In \cite{mohseni2000numerical}, a sharp spectral filter is applied in the azimuthal direction near the axis to remove high-frequency modes and make the code stable. The number of modes that they remove varies with radius from the axis such that the CFL condition is restricted by the radial grid spacing only. 

We consider the 2D diffusion of a scalar in polar coordinates $(r,\theta)$ written in differential form
\begin{equation}
\frac{\partial y}{\partial t} = \frac{1}{r} \frac{\partial}{\partial r} \left( r \frac{\partial y}{\partial r} \right) + \frac{1}{r^2} \frac{\partial^2 y}{\partial \theta^2} \label{eq7_1}
\end{equation}

It is solved in the domain parameterized by $0 \leq r \leq 1$ and $0 \leq \theta < 2\pi$, each direction being uniformly discretized using $N_r$ and $N_\theta$ discretization points respectively. The integral form of the PDE given by Eq.~(\ref{eq7_1}) is solved for each computational element in the form

\begin{equation}
\frac{\partial}{\partial t} \int_{A} y \, dA =  \oint_{d\ell} \frac{\partial y}{\partial n} \, d\ell, \label{eq7_2}
\end{equation}
given by the divergence theorem where $dA$ and $d\ell$ are the respective unit area and unit length respectively. Using second-order finite-volume discretization, where nodal values are stored at the centers of the control volumes, the integrals are approximated using the mid-point rule, which avoids the singularity at $r=0$, and the normal derivative at the cell-face center are approximated using central differencing scheme. Under this discretization, Eq.~(\ref{eq7_1}) becomes

\begin{equation}
\frac{d Y}{d t} = \left( L_\theta + L_r \right) Y + S_r. \label{eq7_3}
\end{equation} 

Consider now the 2D diffusion of the initial solution profile $T(r,\theta,t=0)=\cos{(2\theta)} J_2(\lambda_1 r)$, where $J_2$ is the 2nd-order Bessel function of the first kind and $\lambda_1$ its first root, with known exact solution $T(r,\theta,t)=T(r,\theta,t=0) \exp{(-\lambda_1^2 t)}$. The only required boundary conditions are at $r=1$, and are set to be Dirichlet equal to 0 in accordance with the exact solution.

We use $N_r=10$ and $N_\theta=40$ discretization points in the $r$ and $\theta$ directions, respectively. The resulting maximum time step for stability are $\Delta t_{\mathrm{stability,r}}=6.2 \times 10^{-3}$ and $\Delta t_{\mathrm{stability,\theta}}=4.0 \times 10^{-5}$. We use $\Delta t = 2 \times 10^{-3}$ to reach the final time solution at $t=0.1$ after 50 time steps, which correspond to $\Delta t/\Delta t_{\mathrm{stability,r}}=3.2 \times 10^{-1}$ and $\Delta t/\Delta t_{\mathrm{stability,\theta}}=50$. In order to bypass the numerical stiffness induced by the shrinking length scales near the pole in the azimuthal direction, the TASE methodology can be applied to the operator encoding for diffusion in the $\theta$ direction, resulting in

\begin{equation}
\frac{d Y}{d t} = \left( T_{L_\theta}^{(p)} L_\theta + L_r \right) Y + S_r. \label{eq7_4}
\end{equation} 

Figure~\ref{fig9} shows the qualitative comparison of the solutions obtained with the original explicit scheme with that obtained with the TASE method. The baseline explicit method is fourth-order RK \RevA{(see Table~\ref{But_ERK4} in \ref{But_App})}. While the explicit solution has blown up near the pole, the TASE solution is visually indistinguishable from the exact solution.

\section{Conclusions} \label{conclution}

In this work, a novel TASE method is proposed to handle stiff operators in ODEs and PDEs and to alleviate the inherent strong stability constraint of explicit time integration methods. The proposed method relies on TASE operators that act as preconditioners on the stiff terms, and is flexible to be deployed with most existing explicit time-marching methods. The resulting temporal discretization is still formulated in an explicit manner, yet becomes nearly unconditionally stable. Furthermore, the accuracy order of the baseline explicit methods is preserved as TASE operators can be designed to be arbitrarily high-order accurate. The stiffness from nonlinear convection, diffusion, and source terms can be handled by the corresponding TASE operators within a unified framework, which greatly facilitates the implementation for complex systems of governing equations. Theoretical analyses and stability diagrams show that the $s$-stages $s$th-order explicit RK methods are unconditionally stable when preconditioned by the TASE operators with order $p \leq s$ and $p \leq 2$. On the other hand, the $s$th-order RK methods preconditioned by the TASE operators with order $p \leq s$ and $p > 2$ are nearly unconditionally stable. Based on stability arguments, the single parameter $\alpha$ in TASE operators can be determined based on a prescribed constraint. \RevUs{Contrary to classical implicit methods, the TASE framework allows for solving non-linear problems with arbitrary order without requiring solving a nonlinear system of equations.} A set of stiff physical problems, including those with over-resolved spatial resolution, inhomogeneous boundary conditions, multiple stiff physical process, nonlinear operators, and grid-induced stiffness, demonstrates the performance of the TASE approach. Numerical results reveal that, for all the considered problems, the targeted high-order accuracy is achieved stably for various explicit methods with considerably large time steps.

Although numerically analyzed and verified for specific RK methods, the lack of rigorous  proof of the stability regime for general high-order preconditioned explicit methods represents the main limitation of the present study and will be further investigated in future work. On the other hand, the increased computational cost may overweigh the significance of deploying the TASE framework to very-high-order time marching methods due the required matrix inversions at each explicit marching substep. The reduction of this computational overhead is our ongoing work.

\bibliographystyle{elsarticle-num}


\section*{Acknowledgements}
This investigation was funded by the Advanced Simulation and Computing (ASC) program of the US Department of Energy's National Nuclear Security Administration via the PSAAP-II Center at Stanford (DoE Grant \# 107908). We thank Jack-Michel Cornil for his helpful comments.

\appendix

%
\section{\RevA{Butcher tableaux of Runge-Kutta methods}}
\label{But_App}

\RevA{In this appendix, we report the Butcher tableaux of the Runge-Kutta methods used throughout the paper.}

\begin{table}[!ht]
\begin{center}
\begin{tabular}{ c | c}
0 & 0 \\
\hline
& 1 
\end{tabular}
\end{center}
\caption{Butcher tableau of explicit RK1 (Euler) method used in Section~\ref{Numerical_validations}.} \label{But_ERK1}
\end{table}

\begin{table}[!ht]
\begin{center}
\begin{tabular}{ c | c c}
0 & 0 & 0 \\
1/2 & 1/2 & 0 \\
\hline
& 0 & 1 
\end{tabular}
\end{center}
\caption{Butcher tableau of explicit RK2 (midpoint method) used in Section~\ref{Numerical_validations}.} \label{But_ERK2}
\end{table}

\begin{table}[!ht]
\begin{center}
\begin{tabular}{ c | c c c}
0 & 0 & 0 & 0 \\
1/2 & 1/2 & 0 & 0 \\
3/4 & 0 & 3/4 & 0 \\
\hline
& 2/9 & 1/3 & 4/9 
\end{tabular}
\end{center}
\caption{Butcher tableau of explicit RK3 (Ralston's method) used in Section~\ref{Numerical_validations}.} \label{But_ERK3}
\end{table}

\begin{table}[!ht]
\begin{center}
\begin{tabular}{ c | c c c c}
0 & 0 & 0 & 0 & 0 \\
1/2 & 1/2 & 0 & 0 & 0 \\
1/2 & 0 & 1/2 & 0 & 0 \\
1 & 0 & 0 & 1 & 0\\
\hline
& 1/6 & 1/3 & 1/3 & 1/6 
\end{tabular}
\end{center}
\caption{Butcher tableau of explicit RK4 (``original method") used in Section~\ref{Numerical_validations}.} \label{But_ERK4}
\end{table}

\begin{table}[!ht]
\begin{center}
\begin{tabular}{ c | c}
1 & 1 \\
\hline
& 1 
\end{tabular}
\end{center}
\caption{Butcher tableau of SDIRK1 used in Section~\ref{Numerical_validations}.} \label{But_SDIRK1}
\end{table}

\begin{table}[!ht]
\begin{center}
\begin{tabular}{ c | c c}
0 & 0 & 0 \\
1 & 1/2 & 1/2 \\
\hline
& 1/2 & 1/2 
\end{tabular}
\end{center}
\caption{Butcher tableau of Crank-Nicolson method used in Section~\ref{Numerical_validations}.} \label{But_CN}
\end{table}

\begin{table}[!ht]
\begin{center}
\begin{tabular}{ c | c c}
$\gamma$ & $\gamma$ & 0 \\
$1-\gamma$ & $1-2\gamma$ & $\gamma$ \\
\hline
& $\frac{1}{2}$ & $\frac{1}{2}$
\end{tabular}
\end{center}
\caption{Butcher tableau of SDIRK2 (Table (2) in \cite{pareschi2005implicit}) used in Section~\ref{Numerical_validations}, with $\gamma=1-1/\sqrt{2}$.} \label{But_SDIRK2}
\end{table}

\begin{table}[!ht]
\begin{center}
\begin{tabular}{ c | c c c}
$\gamma$ & $\gamma$ & 0 & 0 \\
$\frac{1+\gamma}{2}$ & $\frac{1-\gamma}{2}$ & $\gamma$ & 0\\
$1$ & $-3\gamma^2/2+4\gamma-1/4$ & $3\gamma^2/2-5\gamma+5/4$ & $\gamma$ \\
\hline
& $-3\gamma^2/2+4\gamma-1/4$ & $3\gamma^2/2-5\gamma+5/4$ & $\gamma$
\end{tabular}
\end{center}
\caption{Butcher tableau of SDIRK3 (Equation 229 in \cite{kennedy2016diagonally}) used in Section~\ref{Numerical_validations}, with $\gamma = 0.43586652150845899941601945$.} \label{But_SDIRK3}
\end{table}

\begin{table}[!ht]
\begin{center}
\begin{tabular}{ c | c c c}
$\gamma$ & $\gamma$ & 0 & 0\\
$1/2$ & $1/2-\gamma$ & $\gamma$ & 0\\
$1-\gamma$ & $2\gamma$ & $1-4\gamma$ & $\gamma$\\
\hline
& $\frac{1}{6(1-2\gamma)^2}$ & $\frac{3(1-2\gamma)^2-1}{3(1-2\gamma)^2}$ & $\frac{1}{6(1-2\gamma)^2}$\\
\end{tabular}
\end{center}
\caption{Butcher tableau of N\o{}rsett's three-stage SDIRK4 used in Section~\ref{Numerical_validations}, with $\gamma=1.06858$.} \label{But_SDIRK4}
\end{table}

\end{document}